\documentclass[11pt,
]{article}

\usepackage{amssymb,amsfonts,amstext,amsmath,amsthm,
latexsym,mathrsfs,amsbsy}
\usepackage{bbold}

\usepackage{marginnote}

\usepackage{mparhack}
\usepackage{esvect}
\input{46.sty}

\vyk{
\usepackage{graphicx}
\usepackage{xcolor}

\usepackage
[%
hypertex
]
{hyperref}
}

\begin{document}

\title
{On Harrington's model in which Separation holds but  
Reduction fails at the 3rd projective level, and on some 
related models of Sami}

\author 
{
Vladimir~Kanovei\thanks{IITP RAS and MIIT,
  Moscow, Russia, \ {\tt kanovei@googlemail.com} --- 
contact author. 
Partial support of   RFFI grant 17-01-00705 acknowledged.
}  
\and
Vassily~Lyubetsky\thanks{IITP RAS,
  Moscow, Russia, \ {\tt lyubetsk@iitp.ru}. 
Partial support of Russian Scientific Fund grant 14-50-00150 
acknowledged. 
}
}

\date 
{\today}

\maketitle

\begin{abstract}
In a handwtitten note of 1975, 
Leo Harrington sketched a construction of a model 
of $\zfc$ (no large cardinals or anything beyond $\ZFC$!) 
in which $\mathbf\Pi^1_3$-Separation holds but 
$\mathbf\Sigma^1_3$-Reduction fails. 
The result has never appeared in a journal or book 
publication except for a few of old references.

MSC 03E15, 03E35
\end{abstract}

\parf{Introduction}

The separation property for a pointclass $K$, or 
simply \rit{\dd KSeparation}, is 
the assertion that any two disjoint sets $X,Y$ 
in $K$ (in the same Polish space) can be separated 
by a set in $K\cap \dop K$, where $\dop K$ is the 
pointclass of complements of sets in $K$. 
The reduction property for a pointclass $K$, or 
simply \rit{\dd KReduction}, is 
the assertion that for any two sets $X,Y$ in $K$ 
(in the same Polish space) there exist 
\rit{disjoint} sets $X'\sq X$, $Y'\sq Y$  
in the same class $K$, such that $X'\cup Y'=X\cup Y$. 

It is known classically from 
studies of 
Luzin \cite{lus:ea,lbook},
Novikov \cite{nov1931,nov1935}, 
Kuratowski \cite{kursep}
that Separation holds for projective classes 
$\fs11$ (analytic sets) and $\fp12$, but fails for 
$\fp11$ (coanalytic sets) and $\fs12$,
while Reduction holds for  
$\fp11$ and $\fs12$, but fails for 
$\fs11$ and $\fp12$, and generally 
\dd KReduction implies \dd{\dop K}Separation
by a simple argument. 

As for the higher projective classes, 
Addison \cite{add2,add1} proved that the axiom 
of constructibility $\rV=\rL$ implies that 
Separation holds for projective classes 
$\fp1n\yd n\ge 3$, but fails for 
$\fs1n\yd n\ge 3$  
while Reduction holds for  
$\fs1n\yd n\ge 3$, but fails for 
$\fp1n\yd n\ge 3$. 
On the other hand, by Martin \cite{martAD}, 
the axiom of projective determinacy $\mathbf{PD}$ 
implies that, similarly to projective level $1$, 
$\fs1n$-Separation and 
$\fp1n$-Reduction hold for all odd numbers $n\ge3$, 
and, similarly to projective level $2$, 
$\fp1n$-Separation and 
$\fs1n$-Reduction hold for all even numbers $n\ge4$. 

Apparently not much is known on Separation and 
Reduction for higher projective classes in generic 
models. 
In a handwtitten note \cite[Part C]{h74} (1975), 
Leo Harrington sketched a construction of a model 
of $\zfc$ in which $\mathbf\Pi^1_3$-Separation holds but 
$\mathbf\Sigma^1_3$-Reduction fails. 
The model was a generic extension of $\rL$ by means of 
the almost-disjoint coding of \cite{jsad}, with no 
reference to determinacy, 
large cardinals or anything beyond $\ZFC$. 
The result has never appeared in a journal or book 
publication except for a few rather old references.\snos
{Hinman \cite[p.\ 230, end of Section V.3]{hin}
communicates a much more general consistency result related 
to the principles of Separation and Reduction, 
absent even in \cite{h74}, citing 
a paper of Harrington entitled 
``\rit{Consistency and independence results in 
descriptive set theory}, to appear in 
\rit{Ann.\ of Math.}'', 
which has apparently never materialized.
Moschovakis \cite[Theorem 5B.3 on p.\ 214]{mDST}
mentions another Harrington's model, 
present in \cite{h74}, in which Separation fails for 
both $\mathbf\Pi^1_3$ and $\mathbf\Sigma^1_3$. 
Another similar reference, to Harrington's models 
in which $\fs1n$-Separation and $\fp1n$-Separation 
both fail for a given $n$, 
see Mathias~\cite[p.\ 166]{matsur}, a comment on P 3110. 
Sami \cite[Thm 1.21]{samiPHD} presents the following 
result with reference to Harrington: 

{\it 
If\/ $n\ge3$, then there exist generic extensions\/ 
$\gN_1\yi\gN_2\yi\gN_3\yi\gN_4$ of\/ $\rL$ such that 
\ben
\aenur 
\itla{sh1}
$\fp1n$-Separation and\/ $\fs1n$-Separation fail 
in\/ $\gN_1\,;$  

\itla{sh2}
$\fp1n$-Separation holds but\/ $\fs1n$-Separation 
fails in\/ $\gN_2\,;$  

\itla{sh3}
$\fp1n$-Separation fails but\/ $\fs1n$-Separation 
holds in\/ $\gN_3\,;$  

\itla{sh3}
$\fp1n$-Separation and\/ $\fs1n$-Separation hold\/ 
in $\gN_4$. 
\een
In addition, there exists a generic extension\/   
$\gN$ of\/ $\rL$ such that
\ben
\aenur 
\atc\atc\atc\atc
\itla{sh5}
$\fp1n$-Separation and $\fs1n$-Separation fail in $\gN$ 
for all\/ $n\ge3$. 
\een
}
\noi
Here, $\gN_1$ is defined in \cite[Part B]{h74} for $n=3$, 
and a hint is given regarding the general case.
A different model, in which 
both $\fs13$-Separation and $\fp13$-Separation fail, 
has recently been defined in \cite{kl28}. 
As for $\gN_2$, the constructible universe itself works 
by Addison. 
Models $\gN_3$ and $\gN_4$ are absent in 
\cite{h74}, generally no generic extensions of $\rL$ 
are known in which $\fs1n$-Separation holds for 
at least one $n\ge 3$.
However a generic model in which both 
(lightface) $\ip1n$-Separation and $\is1n$-Separation 
hold \rit{for sets of integers\/} is given in 
\cite[Part D]{h74}. 
Finally, the existence of a model $\gN$ for \ref{sh5}
is characterized in \cite[Part B]{h74} as an 
\lap{expression of belief}.
} 

Here we present a proof of Harrington's theorem. 
%We prove 

\bte
[Harrington \cite{h74}, Part C]
\lam{mt}
There exists a set-generic extension of\/ $\rL$, 
in which\/ $\fp13$-Separation holds but\/ 
$\fs13$-Reduction fails, and moreover, there is 
a pair of lightfsce\/ $\is13$ sets of reals, 
not reducible 
by a pair of\/ $\fs13$ subsets. 
\ete

In the proof, we'll follow, more or less, 
the flow of Harrington's 
arguments, filling in details and claims wherever 
(we find it) necessary. 
We'll try to preserve even Harrington 
wording wherever possible.
Of most notable deviations, we change Harrington's 
Boolean-valued approach to the poset forcing approach, 
as we observed that the non-absoluteness of the RO 
operation causes problems in understanding of 
the behaviour of certain BAs in different models. 
Of notable additional details, we adjoined some amount 
of definitions and results related to intermediate 
sumbodels of generic extensions, necessary to fully 
understand 
the arguments but near completely avoided (or just hinted) 
in \cite{h74}.

The following is Harrington's comment to Theorem~\ref{mt}  
in \cite[Part C]{h74}.  
\begin{quote}\rit{The above proof was directly inspired by a 
result of Sami, namely: 
there is a model of\/ $\ZFC$ in which\/ 
$\sep(\fp13,\fd13)$ holds but\/ 
$\red(\is13,\is13)$ fails for sets of reals.}
(Note the lightface $\is13$ twice in the second part, 
so $\red(\is13,\is13)$ is $\is13$-Reduction in the above 
sense.) 
\end{quote}
Thus result indeed can be found in Ramez Sami's 
PhD Thesis \cite[Theorem 1.7]{samiPHD}, 
but it has never been published. 

The next theorem presents some related
results in \cite[Thms 1.7,1.18,1.20]{samiPHD}.

\bte
\label{mt'}
\ben
\Renu
\itlb{mt'1}%
It is true in any extension of\/ $\rL$ by a single 
Cohen-generic real that\/ $\is13$-Reduction fails, 
$\ip13$-Separation holds, and if\/ $n\ge4$ then\/ 
$\is1n$-Reduction holds, and hence\/ 
$\fs1n$-Reduction holds as well.~\snos
{To prove that $\is1n$-Reduction implies  
the boldface $\fs1n$-Reduction, it suffices to use  
a double-universal pair of $\is1n$ sets, as those 
used in a typical proof that $\fs1n$-Reduction 
and $\fs1n$-Separation contadict each other.  
This argument does not  work 
for Separation though.}

\itlb{mt'2}%
It is true in any extension of\/ $\rL$ by\/ $\ali$  
Cohen-generic reals that if\/ $n\ge3$ then\/ 
$\is1n$-Reduction holds, and hence\/ 
$\fs1n$-Reduction holds as well.

%\itlb{mt'2}%
%the extension of\/ $\rL$ any number of Solovay-random reals$;$

\itlb{mt'3}%
The same is true in the the Solovay model, \ie, 
the Levy-collapse extension of\/ $\rL$ via an 
inaccessible cardinal.
\een
\ete

We sketch the proof of 
claim \ref{mt'2} in the end of the paper. 
Note that \ref{mt'2} also holds in models obtained 
by adding any uncountable (not necessarily $\ali$) 
number of Cohen-generic reals. 
(Because they are elementarily equivalent to the 
extension by $\ali$ Cohen reals.)
And \ref{mt'2} also holds in extensions by $\ali$ 
or more Solovay-random reals.

\parf{Almost disjoint preliminaries}
\las{nr}

Some definitions related to the almost disjoint 
forcing of \cite{jsad}.
\bit
\item
$\dC=\bse$ (the Cohen forcing). 

\item 
$\La$ 
(the empty string) is the weakest condition in $\dC$. 

\item
$\nse=\ens{s_j}{j<\om}$ is a fixed recursive enumeration.

\item
if $f\in\bn$ then 
$S(f)=\ens{j<\om}{s_j\su f}$.

\item
$\zfcm$ is \ZFC\ without the Power Sets axiom,

%\item
$\bT$ is $\zfcm$ plus $\rV=\rL$ and 
\lap{all sets are countable}.

%\item
%$\btp$ is $\zfcm$ plus 
%\lap{every set belongs to a transitive model of $\bT$}.

\item
$\HC=$ all hereditarily-countable sets.  

\item
$\ang{\xi_\al,n_\al}$ is the $\al$th element of 
the set $\omi\ti\om$, ordered lexicographically.
\eit

\bdf
\lam{fadef}
Reals $f_{\al}\in\bn$ are defined in $\rL$ 
%in 4.5 of \cite{js} 
by induction on $\al<\omi$. 
Suppose that $f_\ga$ are defined for all $\ga<\al$. 
Let $\rL_{\mu(\al)}\mo\bT$ be the smallest model 
containing the sequence $\ga\mto f_\ga$ of 
already defined reals. 
Let $f_\al$ be the Goedel-least real $f\in\bn,$ 
Cohen-generic over $\rL_{\mu(\al)}$ and satisfying 
$s_{n_\al}\su f$.

If $\xi=\xi_\al$ and $n=n_\al$ then let 
$f_{\xi n}:=f_\al$, hence $s_n\su f_{\xi n}$ always holds.
\edf

%Let $\vf=\sis{f_\al}{\al<\omi}$, a sequence of 
%reals $f_\al\in\bn$ in $\rL$.

\vyk{
\ble
\lam{deff}
The sequences\/ $\sis{f_\al}{\al<\omi}$ and\/ 
$\sis{f_{\xi n}}{\xi<\omi,n<\om}$ are\/ $\id{\lomi}1$.
\qed
\ele
}

If $F\sq\bn$ then $\js(F)$ is the corresponding 
almost-disjoint forcing, which consists of all 
pairs $\ang{u,S}$, where $u\sq\om$ and 
$S\sq S(F)=\ens{S(f)}{f\in F}$ are finite sets, 
ordered so that $\ang{t,S}\leq\ang{t',S'}$
(the smaller condition is stronger) 
iff $t'\sq t$, $S'\sq S$, and 
$u\cap A=u'\cap A$ for all $A\in S'$. 

\bit
\item
If $g\in\dn$ then let  
%$\js_g=\js(\ens{f_{\xi i}}{\xi<\omi\land g(i)=0})$.
%$\js_g=\js(F_g)$, where 
$F_g=\ens{f_{\xi i}}{\xi<\omi\yi i<\om\yi g(i)=0}$.

\item
If $e\in\bse$ then let  
$F_e=\ens{f_{\xi i}}{\xi<\omi\yi i<\lh e\yi e(i)=0}$. 
%and $\js_e=\js(F_e)$.

%\item
%We let $\js=\prod_i\js_i$, with finite support.

\item
Define \ $T(g,a)$ \ iff \ 
$\kaz\xi<\omi \:\kaz i\in\om\:
(\text{$S(f_{\xi i})\cap a$ is finite} \leqv g(i)=0)$.
\eit

\ble
[see \cite{jsad}]
\lam{adf}
If\/ $g\in\dn$ in a set universe\/ $\rV$ then 
the forcing\/ $\js(F_g)$ adjoins a real\/ $a\sq\om$ 
satisfying\/ $T(g,a)$.
\qed
\ele 

\bdf
\lam{qqq}
Let $Q\in\rL$ be the forcing notion 
responsible for the following two-step generic 
extension:   
1st, 
we extend a ground set universe $\rV$ 
by a real $g\in\dn$ 
Cohen-generic over $\rV$,  
and 2nd, we extend $\rV[g]$ by $\js(F_g)$.

Thus 
%we may assume that 
$Q$ consists of all triples 
$p=\ang{e,u,S}$, 
where $e\in\bse=\dC$ (a Cohen condition)
while $\ang{u,S}\in \js(F_e)$.
The order is defined so that 
$p=\ang{e,u,S}\leq p'=\ang{e',u',S'}$ 
($p$ is stronger) 
iff $e'\sq e$ and 
$\ang{u,S}\leq \ang{u',S'}$ in $\js(F_e)$.
Note that $1=\ang{\La,\pu,\pu}\in Q$ is the 
largest (and weakest) element in $Q$.

Let $\gQ=Q^\om$ (a finite-support product),
with the product order ${\leq}:={\leq_\gQ}$; 
$p\leq q$ still means that $p$ is stronger. 
Thus $\gQ=\stk\gQ\leq$ is a forcing in $\rL$.
Its largest (= weakest) element $\bon\in\gQ$ 
is defined by $\bon(k)=1$, $\kaz k$.
\edf

\ble
[definability]
\lam{deF}
The sequences\/ $\sis{f_\al}{\al<\omi}$ and\/ 
$\sis{f_{\xi n}}{\xi<\omi,n<\om}$ are\/ $\id{\lomi}1$, 
%over\/ $\lomi$, 
hence\/ $\is{\HC}1$.
The sets\/ $Q$ and\/ $\gQ$ are\/ $\id{\lomi}1$, 
%over\/ $\lomi$, 
hence\/ $\is{\HC}1$.
The relations of compatibility and incompatibility in\/ 
$Q$ and\/ $\gQ$ are\/ $\id{\lomi}1$, 
%over\/ $\lomi$, 
hence\/ $\is{\HC}1$.
\qed
\ele
\bpf 
To circumvent the naturally reqired $\kaz$ over the given 
po set in the definition of incompatibility, define the 
binary operation $\land$ on $Q$ as follows. 
If $p=\ang{e,u,S}$ and $p'=\ang{e',u',S'}$ belong to $Q$ 
then put $p \land q=\ang{e\land e',u\cup u',S\cup S'}$, 
where $e\land e'= e$, or $=e'$, or $=\La$ 
(the empty string) in cases pesp.\ 
$e'\sq e$, $e\sq e'$, 
or $e,e'$ are incomparable in $\dC=\bse.$
Extend $\land$ to $\gQ$ componentwise. 
Then, both in $Q$ and in $\gQ$, conditions 
$p,q$ are incompatible, 
in symbol $p\inc q$, 
iff $p\land q\not\leq p$ or $p\land q\not\leq q$.
This yields the result required.
\epf

\vyk{
\bre
\lam{inco}
If $F\sq\bn$ and $p=\ang{u,S}$, $q=\ang{u',S'}$
belong to $\js(F)$ then we define  
$p\land q=\ang{u\cup u',S\cup S'}$; 
$p\land q\in\js(F)$ as well. 
Note that $p,q$ are compatible in $\js(F)$, 
in symbol $p\com q$, 
iff $p\land q\leq p$ and $p\land q\leq q$. 
Accordingly, there are similar absolute 
operations $\land$ 
on the derived forcing notions $Q$ and $\gQ$, 
such that any conditions $p,q$ are compatible 
in, resp., $Q\yi \gQ$ iff 
$p\land q\leq p$ and $p\land q\leq q$. 
\vyk{
$p\inc q$ (incompatible), 
iff $p\land q\not\leq p$ or $p\land q\not\leq q$.
It follows that
the relations of incompatibility $p\inc q$   
in $Q$, and hence in $\gQ$ as well, are 
$\id{\lomi}1$, hence\/ $\is{\HC}1$, so that 
the expected universal quantifiers over $Q$ and 
$\gQ$ are eliminated. 
}
\ere
}

\bre
\lam{gq}
%$\gQ\in\rL$ and 
(A) 
The forcing $\gQ$ adjoins sequences of the form 
$A=\sis{g_n,a_n}{n<\om}$, 
where each pair $\ang{g_n,a_n}$ is \dd Qgeneric, 
hence $g_n\in\dn$ is Cohen-generic, 
$a_n\sq\om$, and 
$T(g_n,a_n)$ holds.
If $G\sq\gQ$ is generic over a 
set universe $\rV$, and 
$A=A[G]=\sis{g_n[G],a_n[G]}{n<\om}$ is the 
corresponding sequence as above, then $\rV[G]=\rV[A[G]]$. 

(B) 
Any $A$ of such a form can be  
converted to a real $r(A)\sq\om$ by means of any  
recursive bijection between $(\dn\ti\pws\om){}^\om$ 
and $\pws\om$, thus essentially $\gQ$ adds a real,  
so that if $G\sq\gQ$ is generic then 
$r(A[G])\sq\om$ and $\rV[G]=\rV[r(A[G])]$.

(C)
The forcing notions\/ $\js(F_g)$, $Q$, and any 
finite-support product of\/ $Q$, in particular\/ 
$\gQ=Q^\om$ and any $\gQ^\la$, 
satisfy CCC, see \eg\ \cite[Lemma 1 in 4.6]{jsad}. 
\ere

The next lemma is established in \cite{jsad} in 
a somewhat different but pretty similar case,
as a theorem in Section 4.8, pp.\ 95--97, so 
we skip the proof. 

\ble
\lam{c1}
If\/ $g\in\dn$ is Cohen-generic over\/ $\rL$ and\/ 
$B=\sis{\ang{g_n,a_n}}{n<\om}$ is\/ 
\dd\gQ generic over\/ $\rL[g]$ then it holds in\/ 
$\rL[g,B]$ that\/ $\neg\:\sus a\:T(g,a)$.
\qed
\ele
%\bpf[{\cite[Section 4.8]{jsad}}]

\vyk{
Let $A=\sis{a_n}{n<\om}$,  
$G=\sis{g_n}{n<\om}$, and 
$R=\prod_n\js(F_{g_n})\in\rL[G]$.
We make use of the fact that 
$\rL[g,B]=\rL[g,G,A]$, where $\ang{g,G}$ is 
\dd{(\dC\ti\dC^\om)}generic over $\rL$ while 
$A$ is \dd Rgeneric over $\rL[g,G]$. 
Assume towards the contrary that some $r\in R$   
\dd Rforces $T(g,t)$ for some \dd Rname $t\sq R\ti\om$ 
in $\rL[g,G]$.

%We will mimick the argument given in \cite{jsad}.

Let $\rL_\eta$ be a big enough countable elementary 
submodel of $\rL_{\ali}$, such that in particular 
$p,t\in\rL_\eta[g,G]$. 
%such that $\rL_\eta[\ang{g_i}_{i\in\om}]$ 
%contains everything needed 
%to make the following argument work. 
%
Pick $\xi\yi k$ such that $f_{\xi,k}\nin\rL_\eta$ and 
$g(k)=0$, and then a condition $q\in R\yi q\le r$, 
such that 
$q\for_R \text{\lap{$t\cap S(f_{\xi k})$ is finite}}$ 
over $\rL[g,G]$.
We have $q\in R(n_0)=\prod_{n\le n_0}\js(F_{g_n})$
for some $n_0$.

There is a finite set $F$ of reals 
$f_{\da j}$,  
not containing $f_{\xi k}$,  
%does not appear among $\vec f$, 
and such that $F\cap\rL_\eta=\pu$ 
and $q\in\rL_\eta[g,G,F][f_{\xi k}]$.
Note that $f_{\xi k}$ is Cohen-generic over 
$N=\rL_\eta[g,G,F]$.

We may view $q$ as $q(f_{\xi k})$, \ie, $q(f_{\xi k})$ 
is a name in $N$, applied to $f_{\xi k}$. 
%$q$ is in $R$.
%Hence for some $n\ge1$, $q$ is in $\sum_{1\le i\le n}R_i$. 
Hence for all $f_{\da j}$ not in $N$, $q(f_{\da j})$ 
will be in $R(n_0)$ and will extend $p$, 
{\ubf provided} 
$g_i(j)=g_i(k)$, for $i\le n_0$.

Since 
$q\for \text{\lap{$t\cap S(f_{\xi,k})$ is finite}}$, 
and since 
$f_{\xi,k}$ is \dd Cgeneric over $N$, there is $\vec f\in C$ 
such that $\vec f\sq f_{\xi,k}$, and such that 
$$
\vec f\for
\text{\lap{$q(\dof) \for 
(\text{{$t\cap S(\dof)$ is finite}})$%
}}.
$$

Since $\ang{g_i}_{0\le i\le n}$ is 
\dd{\sum_{0\le i\le n}C}generic over $\rL$, 
we can find $f_{\da,j}$ such that $f_{\da,j}\nin N$, 
and $\vec f\sq f_{\da,j}$, and for all $i\yd 1\le i\le n$, 
$g_i(j)=g_i(k)$, and $g_0(j)=1$.

Thus $q(f_{\da,j})$ is in $R$, 
and $q(f_{\da,j})$ extends $p$, and 
$$
q(f_{\da,j}) \for 
\text{\lap{
$t\cap S(f_{\da,j})$ is finite
}}.
$$
Thus $q(f_{\da,j}) \for 
\text{\lap{$\neg\:T(g_0,t)$}}
$.
Contradiction.
}
%\epf

\ble
\lam{homle}
$\gQ$ is homogeneous in the sense that if\/ 
$p,q\in\gQ$ then there is an order automorphism\/ 
$h$ of\/ $\gQ$ such that\/ $p$ and\/ $h(q)$ are 
compatible. 
Therefore
\ben
\renu
\itlb{homle1}%
if\/ $\vpi$ is a formula with names of elements of 
the ground universe\/ $\rV$ as parameters, and some\/ 
$p\in\gQ$ forces\/ $\vpi$ then\/ $\gQ$ 
(\ie, every\/ $q\in \gQ$) forces~$\vpi\,;$

\itlb{homle2}%
if\/ $\vpi(\cdot)$ is a formula with names 
of elements\/ $\rV$ as parameters, 
$t$ is a\/ \dd\gQ name, 
and some\/ $p\in\gQ$ forces\/ $\vpi(t)$, 
then there is another\/ \dd\gQ name\/ $t'$ such 
that $\gQ$ forces\/ $\vpi(t')$ and\/ $p$ forces\/ 
$t=t'$.
\een
\ele
\bpf
As the supports $|p|,|q|\sq\om$ are finite, there is a 
permutation $\pi$ of $\om$ such that the \dd\pi image of 
$|p|$ does not intersect $|q|$. 
Such a $\pi$ induces $h$ as required. 
Claims \ref{homle1} and \ref{homle2} are well-known 
consequences of the homogeneity.
\epf

\vyk{
Recall that a partially ordered set $P=\stk{P}{{\le}}$  
is \rit{separative} if for all $p,q\in P$, 
if $p\not\le q$ then there exists an $r\le p$ 
that is incompatible with $q$, \cite[14.8]{jechmill}.

\ble
\lam{sepf}
The forcing notions\/ $\dC$ (Cohen), $\js(F_g)$, $Q$, $\gQ$ 
are separative. 
\ele
\bpf
The result for $\js(F_g)$ follows from the fact that if 
$\al,\ga_0,\dots,\ga_n<\omi$ and 
$\al\nin\ans{\ga_0,\dots,\ga_n}$ then the set
$S(f_\al)\bez\bigcup_{k\le n}S(f_{\ga_k})$
is infinite.
\epf
}

\parf{On intermediate models}
\las{inte}

Given a forcing notion $P=\stk{P}{\le}$ 
in a ground set universe $\rV$, 
if a set $X\sq\rV$ belongs to 
a \dd Pgeneric extension $\rV[G]$ of $\rV$, then the 
submodel $\rV[X]\sq\rV[G]$ is a generic 
extension of $\rV$. 
(But it is not asserted that the set 
$X$ itself is generic over $\rV$!) 
This issue has been exhaustively studied in terms of 
boolean-valued forcing (see \eg\ Lemma 69 in \cite{jechE}), 
which we avoid here. 
%due to the non-absuluteness 
%of the operation $\ro P$. 
Instead we make use of the classical \dd\Sg construction 
by Solovay \cite{sol}, rendered here only for 
the case $P=\gQ$ and $X\sq\rV$. 
(See \cite[Section 1]{zapt} or \cite{kl21} for the 
treatment in the case when $X$ is not a subset of the 
ground set universe $\rV$.)
Basically, the results below hold for any $P\in\rL$, and 
if $P\nin \rL$ then the results also hold with $P$ as 
a uniform parameter.

\bdf
[Solovay \cite{sol}]
\lam{ssig}
Assume that $t\in\rV$, $t\sq \gQ\ti\rV$. 
($\rV$ being a ground set universe.) 
Let $X\sq \rV$ be a set in a generic extension  
of $\rV$. 
We define $\sid Xt =\gQ\bez\bigcup_{\al<\vt}W_\al$, 
where the increasing sequence 
of sets $W_\al\sq\gQ$ is defined in $\rV[X]$ 
by induction,  
and an ordinal $\vt$ is determined in the course 
of construction.  
\ben
\nenu
\itlb{ssig1}%
$W_0$ consists of all conditions $p\in\gQ$ such that 
{\ubf either} there is a set $x\in X$ such that $p$ is 
incompatible in ${\gQ}$ with any condition $q$ satisfying 
$\ang{q,x}\in t$, 
{\ubf or} there is $x\nin X$ and 
a weaker condition $q\geq p$ such that $\ang{q,x}\in t$.

\itlb{ssig2}%  
$W_{\al+1}$ consists of all conditions $p\in {\gQ}$ such that 
there is a dense set $D\in\rV$, $D\sq {\gQ}$ satisfying: 
if $q\in D$ and $q\leq p$ then $q\in W_\al$. 

\itlb{ssig3}%  
If $\la$ is limit then $W_\la=\bigcup_{\al<\la}W_\al$.
\een 
Note that $\al<\ba\imp W_\al\sq W_\ba$, hence 
there is an ordinal $\vt$ satisfying $W_\vt=W_{\vt+1}$, 
and then $W_\xi=W_\vt$ for all $\xi>\vt$. 
Finally, let $W=\bigcup_{\al<\vt}W_\al$. 
The set $\sid Xt ={\gQ}\bez W$ 
contains all conditions $p$ which, roughly speaking, 
are compatible with the assumption that $X=t[G]$ 
for a \dd {\gQ}generic set $G$ containing $p$.
\edf

A set $D\sq {\gQ}$ is \rit{dense} iff 
$\kaz p\in {\gQ}\,\sus q\in D\,(q\leq p)$, and 
\rit{open} dense if in addition 
$(q\in D\land q\leq p)\imp p\in D$.
The set 
$t[G]=\ens{x}{\sus p\in G\,(\ang{p,x}\in t)}$ 
is the \rit{\dd{G}valuation} of $t$. 
The next lemma evaluates the length $\vt$ of the 
construction of \ref{ssig}.

\ble
\lam{evt} 
Under the assumptions of\/ \ref{ssig}, $\vt\le\omi$.
\ele
\bpf
To prove $W_{\omi+1}=W_{\omi}$, 
let $p\in {\gQ}$, hence there is a dense set    
$D\in\rV$, $D\sq {\gQ}$, satisfying   
$({q\in D}\land {q\leq p})\imp{q\in W_{\omi}}$. 
%We can assume 
%
%
The set
$$
D'=\ens{q'\in D}
{%\sus r\in A\,(q\leq r) \land
%\big(
q'\leq p\lor q',p\text{ are incompatible}
%\big)
}
$$
is still dense and satisfies $(\ast)$ 
$({q'\in D'}\land {q'\leq p})\imp{q'\in W_{\omi}}$. 
\vyk{
Furthermore the set 
$$
D''=\ens{q''\in\gQ}
{\sus q'\in D'\,(q''\leq q')}
$$
is even \rit{open} dense, and we still have 
$({q''\in D''}\land {q''\leq p})\imp{q''\in W_{\omi}}$. 
(Indeed if $q''\leq q'\in D'$ and $q''\leq p$ then 
$q',p$ are compatible, hence $q'\leq p$, so 
$q'\in W_{\omi}$. 
But, by induction, each set $W_\al$ satisfies
$q''\leq q'\in W_\al\imp q''\in W_\al$.) 
}%

Let $A\sq D'$ be a maximal antichain; 
$A$ is countable by \ref{gq}(C), 
hence there is $\al<\omi$ such that 
$(\dag)$
$A\bigcap W=A\cap W_\al$.
By the maximality of $A$,  the set 
$$
D''=\ens{q''\in\gQ}
{\sus r\in A\,(q''\leq r)}
$$
is dense, and   
even \rit{open} dense.
We claim that 
$(\ddag)$
$({q''\in D''}\land {q''\leq p})\imp{q''\in W_{\al}}$; 
this implies $p\in W_\al\sq W_{\omi}$, ending the 
proof of the lemma. 
Thus prove $(\ddag)$. 

By definition, $q''\leq r$ for some $r\in A$. 
Thus $r,p$ are compatible. 
Therefore, as $A\sq D'$, we have $r\leq p$. 
To conclude, $r\leq p$ and $r\in D'$. 
It follows that $r\in W_{\omi}$ by $(\ast)$, 
hence $r\in W_\al$ by $(\dag)$. 
It follows that $q''\in W_\al$. 
(Indeed, by induction, each set $W_\al$ satisfies
$q''\leq r\in W_\al\imp q''\in W_\al$.)
As required.
\epf

\bte
[Solovay \cite{sol}]
\lam{32s}
Under the assumptions of\/ \ref{ssig}, 
suppose that\/ 
%${\gQ}=\stk{{\gQ}}{\le}$ is a forcing    
%in a set universe $\rV$, $t\in\rV$, $t\sq {\gQ}\ti\rV$,  
a set\/ $G\sq {\gQ}$ is\/ \dd {\gQ}generic over\/ $\rV$, 
and\/ $X=t[G]$. 
Then\/ 
\ben
\renu
\itlb{32s1}%
$G\sq \Sg=\sid Xt$ and\/ $X=t[\Sg]$ ---  
hence\/ $\rV[\Sg]=\rV[X]$, 

\itlb{32s2}%
$G$ is\/ \dd\Sg generic over the intermediate 
model\/ $\rV[\Sg]=\rV[X]\sq\rV[G]$ ---  
hence\/ $\rV[G]$ is a set generic extension of\/ 
$\rV[X]\,;$ 

\itlb{32s3}%
if\/ $G'\sq\Sg$ is\/ \dd\Sg generic 
over\/ $\rV[\Sg]=\rV[X]$ then\/ $G'$ is\/  
\dd {\gQ}generic over\/ $\rV$ and $t[G']=X$. 
% --- it follows that\/ $X=t[\Sg]$.
\qed 
\een
\vyk{
It follows that\/ $\Sg=\sid Xt\rV$ is equal to the set 
of all conditions\/ $p\in {\gQ}$ such that, 
in a suitable collapse-generic extension of\/ $\rV$, 
there is 
a set\/ $G'\sq {\gQ}$, \dd {\gQ}generic over\/ $\rV$ and 
satisfying\/ $p\in G'$ and\/ $t[G']=X$.
} 
\ete

\bcor
\lam{32c}
Under the assumptions of Theorem~\ref{32s}, let\/ 
$\vpi(\cdot)$ be a formula with\/ \dd {\gQ}names 
for sets in\/ $\rV$ allowed as parameters. 
Then\/ $\vpi(X)$ is true in\/ $\rV[X]$ iff there 
is a condition\/ $p\in\Sg$ that\/ \dd {\gQ}forces\/ 
$\rV[t[\uG]]\mo\vpi(t[\uG])$ over\/ $\rV$.
\ecor
\bpf
To prove $\mpi$, the less trivial direction,
assume that a condition $p\in\Sg$ 
\dd {\gQ}forces $\rV[t[\uG]]\mo\vpi(t[\uG])$ over $\rV$.
Consider a set $G'\sq\Sg$, \dd\Sg generic over $\rV[X]$ 
and containing $p$.
Then $G'$ is \dd {\gQ}generic over $\rV$ and $t[G']=X$.
It follows that $\vpi(X)$ is true in\/ $\rV[X]$ 
by the choice of $p$. 
\epf

\bdf
[see \cite{kl32e}]
\lam{kkan}
Let 
%$P=\stk{P}{\le}$ is a forcing notion 
%in a set universe $\rV$, 
$t\in\rV$, $t\sq {\gQ}\ti\rV$. 
Define, in $\rV$, the order relation $\tle t$ 
on ${\gQ}$ as follows: 
$p\tle t q$, iff $p$ \dd{\gQ}forces over $\rV$ that 
$\namy q\in\sid{\namy t[\uG]}{\namy t}$. 

Let 
$\gQ_t=\stk\gQ{\tle t}$.
\edf

Under the hypothesis that for any $p\in P$ there is 
a \dd Pgeneric set $G\sq P$ containing $p$, the 
relation $p\tle t q$ is equivalent to the following: 
if $G$ is a set \dd\gQ generic over $\rV$ and 
containing $p$ then $q\in\sid Xt$. 
\vyk{
Roughly speaking, $p\tle tq$ means that  
every sentence $\vpi(t[\uG])$, 
\dd Pforced by $q$ to be true in 
$\rV[t[\uG]]$, is \dd Pforced by $p$ as well.
}%

The next theorem contains the main application of 
the orders $\tle t$.

\bte
%[proved in \cite{kl32e}]
\lam{32t}
Suppose that\/ 
%$P=\stk{P}{\le}$ is a forcing notion  
%in a set universe $\rV$, and\/ 
$t\in\rV$, $t\sq\gQ\ti\rV$. 
Then\/ 
\ben
\renu
\itlb{32t1}%
${\tle t}$ is a partial order relation 
on\/ $\gQ$ satisfying\/ ${p\le q}\imp {p \tle t q}$,  

\itlb{32t2}%
if a set\/ $G\sq\gQ$ is\/ \dd\gQ generic over\/ $\rV$ 
and\/ $X=t[G]$, then the set\/ $\sid Xt$ is\/ 
\dd{\gQ_t}generic over\/ $\rV$ and\/ $\rV[X]=\rV[\sid Xt]$,  

\itlb{32t3}%
if a set\/ $\Sg\sq\gQ$ is\/ \dd{\gQ_t}generic over\/ $\rV$, 
and a set\/ $G'\sq\Sg$ is\/ 
\dd{\Sg}generic over\/ $\rV[\Sg]$, 
then\/ $G'$ is\/ \dd{\gQ}generic over\/ $\rV$ and\/ 
$\Sg\in\rV[G']$. 
\een
\ete
\bpf
Claims \ref{32t1}, \ref{32t2} are established in 
\cite{kl32e}. 
To prove \ref{32t3}, suppose towards the contrary that  
$p\in\gQ$ forces the negation, that is, if 
$\Sg\sq\gQ$ is \dd{\gQ_t}generic over $\rV$,  
$G'\sq\Sg$ is \dd{\Sg}generic over $\rV[\Sg]$, 
and $p\in G'$, then either $G'$ 
is \text{\ubf not} \dd{\gQ}generic over $\rV$ 
or $\Sg\in\rV[G']$.
%To get a contradiction, c
Consider a set $G'\sq\gQ$, 
containing $p$ and \dd\gQ generic over $\rV$. 
Let $X'=t[G']$.
The set $\Sg'=\sid {X'}t$ is 
\dd{\gQ_t}generic over $\rV$ by \ref{32t2}, 
$G'$ is \dd{\Sg}generic over $\rV[\Sg']$ 
by Theorem~\ref{32s}\ref{32s2}, 
and $X'=t[G']\in\rV[G']$, hence 
$\Sg'=\sid {X'}t\in\rV[G']$ as well.
Finally, $p\in G'\sq\Sg'$, 
which contradicts the choice of $p$. 
\epf

\ble
\lam{sust}
Let\/ $G\sq\gQ$ be\/ \dd\gQ generic over\/ $\rV$, and
$x\sq\om$ be a real in\/ $\rV[G]$. 
Then there is a countable set\/ $t\in\rV$, $t\sq\gQ\ti\om$ 
such that\/ $x=t[G]$ and\/ $\rV[x]$ is a\/ 
\dd{\gQ_t}generic extension 
of\/ $\rV$.
{\rm(But $x$ itself is not asserted to be a generic.)}
\ele
\bpf
By basic forcing theory, there is a set $t\in\rV$, 
$t\sq\gQ\ti\om$, satisfying $x=t[G]$, and by CCC 
(see \ref{gq}(C)) there is a countable such $t$.
Now use Theorem~\ref{32t}\ref{32t2}.
\epf

\bre
\lam{ccc'}
Under the assumptions of the theorem, 
%if a forcing notion 
%$P=\stk{P}{\le}$ satisfies CCC (in $\rV$) 
%then so does 
$\gQ_t$ satisfies CCC (in $\rV$). 
Indeed, as ${(\leq)}\sq{(\tle t)}$, any \dd{\tle t}antichain 
is a \dd{\leq}antichain as well.
\ere

\vyk{
\bre
\lam{t=G}
In the assumptions of the theorem, 
if $P$ forces $\rV[t[\uG]]=\rV[\uG]$, then it can be 
assumed that $\tle t$ coincides with $\le$.
\ere
}

\vyk{
\bdf
\lam{haf}
We return to the forcing notion $\gQ=Q^\om\in\rL$.  
Let $\le$ be the natural partial order 
of $\gQ$ in $\rL$, so that $\gQ=\stk\gQ{\le}$.  
If $t\sq \gQ\ti\om$ 
(a \dd\gQ name for a subset of $\om$), 
%then let $\tle t$ denote $\tlf\gQ {\,t}$, and let 
$\gQ_t=\stk\gQ{\tle t}$.

The \rit{Harrington fan} $\fan\gQ$ consists of 
all forcing notions of the form 
$\gQ_t=\stk\gQ{\tle t}$, 
where $t\sq \gQ\ti\om$ is at most countable. 
\edf
}

\bre
\lam{qqp}
(A)
Let the \rit{Harrington fan} $\fan\gQ$ consist of 
all forcing notions of the form 
$\gQ_t=\stk\gQ{\tle t}$, 
where $t\sq \gQ\ti\om$ is at most countable. 

(B)
Coming back to Remark~\ref{gq}, let $\bta\in\rL$ be a 
canonical \dd\gQ name of the real $r(A[G])$ 
as in \ref{gq}(B), so that if $G$ is \dd\gQ generic 
then $\bta[G]=r(A[G])$ and hence $\rV[G]=\rV[\bta[G]]$, 
and moreover if $G\ne G'$ then $\bta[G]\ne \bta[G']$. 
This implies that the order $\tle\bta$ on $\gQ$ 
coincides with $\leq$, but by means of a rather legthy 
argument, which includes the verification of the 
separativity of the forcing notion $\gQ$. 
In order to circumwent these complications, 
{\ubf\boldmath 
it will be outright assumed that the partial order 
$\tle\bta$ on $\gQ$ coincides with $\leq$}, 
for this particular \dd\gQ name $\bta$, 
and accordingly if a set $G\sq\gQ$ is generic and 
$\bta[G]=r(A[G])=\br\sq\om$ then 
$\sid\br\bta=G$, so \ref{32s}\ref{32s1} still holds. 
With this amendment, we have $\gQ\in\fan\gQ$.\vom

(C)
In the notation of Remark~\ref{gq}(A), 
let $\bc\in\rL$ be the 
canonical \dd\gQ name for the set $\ens{k<\om}{g_0(k)=0}$, 
so that $\gQ$ forces that $\bc[\uG]\sq\om$ 
is Cohen-generic over the 
ground universe. 
Then the forcing $\gQ_\bc\in\fan\gQ$ adds a Cohen real. 
\ere

\parf{Absoluteness of the $\Sg$ construction}
\las{absol}

We have to consider a subtle issue related to 
the construction of $\sid Xt$, namely, 
its formal dependence 
of the choice of $\rV$ in (2) of Definition~\ref{ssig}. 
The next lemma shows that the dependence 
can be eliminated in a really important case.

\ble
\lam{abso*}%
Under the assumptions of Definition~\ref{ssig}, 
suppose that, in addition, 
$\rV$ is a set-generic extension of\/ $\rL[t]$.
Then\/ $\sie Xt \rV=\sie Xt{\rL[t]}$.  
\ele
\bpf
Assume that $\Pi\in\rL[t]$ is a forcing notion, 
and $\rV=\rL[t][H]$, 
where $H$ is \dd\Pi generic over $\rL[t]$. 
Then $H$ is \dd\Pi generic over $\rL[t][X]$ by 
the product forcing theorem.
Prove by induction that $W^\rV_\al=W^{\rL[t]}_\al$. 

It suffices to handle the inductive step 
$\al\to\al+1$ in \ref{ssig}(2). 
Thus suppose that $W^\rV_\al=W^{\rL[t]}_\al=W_\al$ and 
prove $W^\rV_\alp=W^{\rL[t]}_\alp$.
As $\rL[t]\sq\rV$, we have 
$W^{\rL[t]}_\alp\sq W^\rV_\alp$.
To prove the opposite inclusion, suppose that 
$p_0\in W^{\rV}_\alp$, and this is witnessed by a 
dense set $D\in\rV$, $D\sq {\gQ}$, as in \ref{ssig}(2). 
The goal is to prove $p_0\in W^{\rL[t]}_\alp$.
We have $D=\tau[H]$, where $\tau\in \rL[t]$, 
$\tau\sq \Pi\ti {\gQ}$ 
(a \dd\Pi name of a subset of ${\gQ}$).  
There is a condition $\pi_0\in H$, which \dd\Pi forces, 
over $\rL[t][X]$, that 
%\pagebreak[0] 
$$
\text{\lap{$\tau[\uH]$ is dense and 
$\kaz p\in\tau[\uH]\,(p\leq p_0\imp p\in W_\al^{\smile})$
},}
\eqno(\dag)  
$$
where $W_\al^{\smile}=\Pi\ti W_\al$ is the canonical 
\dd\Pi name for the set $W_\al\in \rL[t][X]$.

We can wlog assume that $(\ddag)$ 
$\ang{\pi,p}\in \tau\land \pi'\in\Pi\land 
\pi'\le_{\Pi} \pi\imp \ang{\pi',p}\in \tau$.

We notice that the sets 
$D'_1=\ens{p\in {\gQ}}
{p\leq p_0\land\sus\pi\in\Pi\,(\pi\le_\Pi\pi_0
\land\ang{\pi,p}\in\tau)}$, 
$D'_2=\ens{p\in {\gQ}}{p_0,p\text{ are incompatible}}$, 
and $D'=D'_1\cup D'_2$  
belong to $\rL[t]$.

We claim that $D'$ is dense in ${\gQ}$. 
Indeed let $p\in {\gQ}$. 
If $p$ is incompatible with $p_0$ then immediately $p\in D'$. 
If otherwise, then we can assume that $p\leq p_0$. 
As $\pi_0$ forces $(\dag)$, there is a condition $\pi\in H$, 
$\pi\le_\Pi\pi_0$, and some $p'\in {\gQ}$, $p'\leq p$, 
such that 
$\pi$ forces $\pi'\in\tau[\uH]$ --- that is, 
$\ang{\pi,p'}\in\tau$ by $(\ddag)$. 
Then $p'\in D'$, as required. 

We finally claim that if $p\in D'$ and $p\leq p_0$ then 
$p\in W_\al$. 
Indeed, $p\nin D'_2$, hence, $p\in D'_1$. 
Let this be witnessed by $\pi\in\Pi$, $\pi\le_\Pi \pi_0$. 
Then $\pi$ obviously forces $p\in\tau[\uH]$, and hence, 
as $\pi$ also forces $(\dag)$, we conclude that 
$p\in W_\al$, as required.

Thus $D'$ witnesses that $p_0\in W^{\rL[t]}_\alp$, 
as required.
\epf

\bcor
\lam{absc}%
Under the assumptions of Lemma~\ref{abso*}, if\/ 
$\rV_1=\rL[t,Y_1]$ and\/ $\rV_2=\rL[t,Y_2]$,  
where\/ $Y_1\yi Y_2\in\rV$ and\/ 
$Y_1\cup Y_2\sq\rL[t]$, then\/ 
\ben
\renu
\itlb{absc1}%
$\sie Xt{\rV_1}=\sie Xt{\rV_2}=
\sie Xt{\rV}=\sie Xt{\rL[t]}$.

\itlb{absc2}%
${\tle t^{\rV_1}} = {\tle t^{\rV_2}} =
{\tle t^{\rV}} = {\tle t^{\rL[t]}}$.
\een  
\ecor
\bpf
%The additional assumption of the lemma is not a burden 
%since 
If $\rV$ is a set-generic extension 
of $\rL[t]$ then any subextension $\rL[t,Y]$, 
where $Y\in\rV\yi Y\sq\rL[t]$, 
is a set-generic extension 
of $\rL[t]$ as well by \cite{gri} or Theorem~\ref{32t}. 
This implies \ref{absc1} by Lemma~\ref{abso*}, 
and then \ref{absc2} also follows by a routine argument.
%$\rV$ is a set-generic extension of 
%any $\rL[Y]$, where $Y\in\rV$ and $Y\sq\rL$, by 
%Theorem~\ref{32s}\ref{32s3}.
\epf

\bbla
\lam{blan}
We'll freely 
use the notation $\sid Xt$ and $\tle t$ without 
reference to the ground universe,  
due to Corollary \ref{absc}.
Indeed, the universes considered will always 
be subuniverses of a fixed generic extension 
of $\rL$.
\ebla

\bcor
\lam{ptsr}%
\sloppy
Under the assumptions of Lemma~\ref{abso*}, 
the relation\/ 
$\tle t$ belongs to\/ $\rL[t]$, 
and we have\/  
$\rL[t,X]=\rL[t,\sid Xt]$.\qed
\ecor

\parf{Definability of the $\Sg$ construction}
\las{defsol}

%Recall that $\HC=$ all hereditarily-countable sets.
%
Consider the sets 
${\bE}=\ens{A\in\rL}{
%A\in\rL\land 
A\text{ is a maximal antichain in }\gQ}$
and 
$$
\bbsg
=
\ens{\ang{X,t,p}}
{X\sq\om\land t\in\lomi\yi t\sq\gQ\ti\om %\text{ countable}
\land{p\in\sid Xt}}. 
$$
We have $\bbsg\sq\lomi$, and also 
${\bE}\sq\lomi$ since $\gQ$ is CCC.

\ble
\lam{Wdef}
%$W\sq\lomi$ and\/ 
${\bE}$ is\/ 
%$\is{}1$ over\/ $\rL_{\omi}$ hence\/ 
$\is{\lomi}1$.
\ele
\bpf
See Section~\ref{nr} on $\zfcm$ and $\bT$.
Let $\btp$ be the theory $\zfcm$ 
%(minus the Power Sets) 
plus the axiom saying that every set belongs to 
a countable model $\rL_\al\mo\bT$. 
We claim that the following are equivalent: 
(1) $A\in {\bE}$,\vom 

(2) 
$\sus\la<\omi(\rL_\la\mo\btp\land A\in\rL_\la\land 
A\text{ is a maximal antichain in }\gQ\cap\rL_\la)$,\vom 

(3)
$\kaz\la<\omi(\rL_\la\mo\btp\land A\in\rL_\la\imp 
A\text{ is a maximal antichain in }\gQ\cap\rL_\la)$.\vom

\noi
If this is established then Lemma~\ref{deF}  leads 
to the definability result required.

It remains to prove the claim. 
$(1)\imp(3)$ is obvious.  
%the nontrivial implication. 
To prove 
$(3)\imp(2)$, 
consider any countable elementary submodel 
$\rL_\la$ of $\rL_{\omi}$ containing $A$. 
%Let $f:M \onto \rL_\la$ be the Mostowski collapse. 
Then $\rL_\la$ is a model of $\btp$ (as so is $\rL_{\omi}$). 
%and $f(\gQ)=\gQ\cap\rL_\la$. 
Thus $A$ is a maximal antichain in $\gQ\cap\rL_\la$.

Now prove $(2)\imp(1)$, the nontrivial implication. 
Suppose that $\la<\omi$, $A\in\rL_\la\mo\btp$, 
and $A$ is a maximal antichain in $\gQ\cap\rL_\la$. 
Then $A$ is an antichain in $\gQ$, since being antichain 
means that $p\land q\not\leq p$ or $\not\leq q$ for any 
$p\ne q$ in $A$, and $\land$ is an absolute operation  
(see the proof of Lemma~\ref{deF}). 
It remains to prove that $A$ is a maximal antichain. 
\vyk{
This takes some space. 
First of all, let $\nu=\omi^{\rL_\la}$. 
Then in fact $A\in\rL_{\nu}$, 
}

Suppose to the contrary 
that some $p\in\gQ\bez\rL_\la$ is incompatible with 
every $q\in A$. 
By definition there are finitely many  
reals of the form $f_{\xi n}$ with $\xi\ge\la$, 
occurring in $p$, so we may write 
$p=p(f_{\xi_0, n_0},\dots,f_{\xi_k, n_k})$, $k<\om$,  
$\xi_i\ge\la$ for all $i$. 
(The substitution form 
$p(f_{\xi_0, n_0},\dots,f_{\xi_k, n_k})$
is naturally chosen so that 
$p(f_{\et_0, n_0},\dots,f_{\et_k, n_k})\in\gQ$ 
for any other string of ordinals 
$\et_0,\dots,\et_k<\omi$.) 
The reals $f_{\xi n}$ with $\xi<\la$ may occur as well, 
but they belong to $\rL_\la$ and are not  
to be explicitly mentioned. 
But anyway there is an ordinal $\nu<\la$ such that 
all $f_{\xi n}$ with $\xi<\la$, occuring in $p$, actually 
satisfy $\xi<\nu$, and in addition $A\in\rL_\nu\mo\bT$. 

Note that by construction the string of reals 
$\ang{f_{\xi_0, n_0},\dots,f_{\xi_k, n_k}}\in(\bn){}^{n+1}$ 
is Cohen generic over $\rL_\la$, hence over $\rL_\nu$ as 
well. 
Therefore the property 
\bce
\lap{$p(f_{\xi_0, n_0},\dots,f_{\xi_k, n_k})$ 
is incompatible with every $q\in A$}
\ece
is forced over $\rL_\nu$, in the sense that there exist 
strings $e_0,\dots,e_k\in\nse,$ such that 
$e_i\su f_{\xi_i, n_i}$, $\kaz i$, 
and if $\ang{y_0,\dots,y_n}$ is Cohen generic 
over $\rL_\nu$ with $e_i\su y_i$, $\kaz i$, then still 
$p(y_0,\dots,y_k)$ is incompatible with every $q\in A$. 

It remains to note that, since $\nu<\la$, there exist 
intermediate ordinals $\et_0,\dots,\et_k\in\la\bez\nu$ 
such that the reals $y_i=f_{\et_i, n_i}$ satisfy 
$e_i\su y_i$, $\kaz i$. 
Then $\ang{y_0,\dots,y_n}$ is Cohen generic 
over $\rL_\nu$, hence $p'=p(y_0,\dots,y_k)$ is 
incompatible with each $q\in A$ by the above. 
And on the other hand $p'\in\gQ\cap\rL_\la$, 
a contradiction. 
\epf

%\np

\ble
\lam{Sdef}
%$W\sq\lomi$ and\/ 
$\bbsg$ is definable in\/ 
$\hc=\text{\rm hereditarily-countable sets}$ 
by a conjunction of the\/ $\is{}1$ formula\/ 
{\rm\lap{$t\in\lomi$}} and a\/
$\ip{}1$ formula\/ $\jsg(X,t,p)$.
\ele
\bpf
Assume that $X\sq\om$, $t\in\rL$, $t\sq\gQ\ti\om$, 
and $p\in\gQ$. 
Let a \rit{maximal \dd pantichain} 
be any maximal antichain $A\sq\gQ$ such that if $q\in A$ 
then either $q\yi p$ are incompatible or $q\leq p$.
%Let $t\in\rL\yi t\sq\gQ\ti\omi$.
Come back to the sets $W_\al$ in \ref{ssig}. 
As $t\in\rL$, Lemma~\ref{abso*} allows us to 
consider only dense sets $D\in\rL$ in \ref{ssig}\ref{ssig2}.

If $\al<\omi$ then let an \rit{\dd\al ladder} be any 
sequence $\sis{W'_\xi,\cA_\xi}{\xi\le\al}$ such that 
each $W'_\xi\sq\gQ$ is at most countable, each 
$\cA_\xi\sq\bE$ is at most countable, and   
\ben
\nenup
\itlb{psig1}\msur%
$W'_0\sq W_0$ \ 
(the latter defined as in \ref{ssig}\ref{ssig1});

\itlb{psig2}%
if $\xi+1\le\al$ and $p\in W'_{\xi+1}$ then there is 
a maximal \dd pantichain $A\in\cA_\xi$ such that 
$r\in W'_\xi$ holds for all $r\in A\yi r\leq p$;

\itlb{psig3}%
if $\xi$ is limit then $W'_\xi=\bigcup_{\et<\xi} W'_\et$.
\een
We assert that if $p\in\gQ$ and $\al<\omi$ then: 
\ben
\fenu
\itlb{psig*}
\rit{$p\in W_\al$ 
iff there is an \dd{\al}ladder\/ 
$\sis{W'_\xi,\cA_\xi}{\xi\le\al}$ 
such that\/ $p\in W'_\al$}. 
\een
If this is established then 
$\bigcup_{\al<\omi}W_\al$ becomes 
a $\is\hc1$ set. 
(The incompatibility in \ref{ssig}\ref{ssig1}, to 
which \ref{psig1} refers, is handled by Lemma~\ref{deF}.)
Then $\sid Xt=\gQ\bez \bigcup_{\al<\omi}W_\al$ becomes 
a $\ip\hc1$ set, and the lemma easily follows. 
Note that the union needn't exceed 
$\omi$ by Lemma~\ref{evt}.

In the direction $\mpi$ of \ref{psig*}, 
we prove by induction that 
$W'_\xi\sq W_\xi$. 
The nontrivial step is \ref{psig2}. 
Let $p\in W'_{\xi+1}$, and let this be witnessed by 
$A\in\cA_\xi$ in the sense of \ref{psig2}. 
As $A\in\bE$, the set 
$D=\ens{q\in\gQ}{\sus r\in A\,(q\leq r)}$ is dense 
and $D\in\rL$. 
It remains to prove that if $q\in D\yi q\leq p$, 
then $q\in W_\xi$, see \ref{ssig}\ref{ssig2}. 
Indeed, by construction there is $r\in A$ with $q\leq r$. 
But $A$ is a \dd pantichain, hence 
either $r\yi p$ are incompatible or $r\leq p$. 
However $q\leq p\yd q\leq r$, excluding the `either' case. 
Thus $r\leq p$. 
It follows by the choice of $A$ that $r\in W'_\xi$. 
Thus $r\in W_\xi$ by the inductive hypothesis. 
We conclude that $q\in W_\xi$ as well, since $q\leq p$. 

We prove $\imp$ in \ref{psig*} by induction on $\al$. 
The nontrivial step is still \ref{psig2}. 
Thus suppose that $p\in W_{\al+1}$, and let 
this be witnessed by a dense set 
$D\in\rL\yd D\sq\gQ$ 
in the sense of \ref{ssig}\ref{ssig2}. 
Let $D'$ consist of all $q\in D$ such that 
$q\leq p$ or $q$ is incompatible with $p$; 
then $D'\in\rL$ is still dense and witnesses 
$p\in W_{\al+1}$. 
Consider a maximal antichain $A\sq D'$ in 
$\rL$, so that $A\in\bE$. 
Then 
$A$ is a maximal \dd pantichain by the 
definition of $D'$, and if $q\in A\yd q\leq p$ 
then $q\in W_\al$ by the choice of $D$, hence, 
by the inductive hypothesis, there is  
an \dd{\al}ladder
$\sis{W_\xi(q),\cA_\xi(q)}{\xi\le\al}$ satisfying 
$q\in W_\xi(q)$. 
To accomplish the proof of $\imp$ in \ref{psig*}, 
define an \dd{(\al+1)}ladder by 
$W'_\xi=\bigcup_{q\in A\yi q\leq p}W_\xi(q)$ and 
$\cA_\xi=\bigcup_{q\in A\yi q\leq p}\cA_\xi(q)$ 
for $\xi\le\al$, and separately 
$W'_{\al+1}=\ans p$ and $\cA_{\al+1}=\ans A$. 
\vyk{
Finally, the relation $p\in\sid Xt$ is equivalent 
to $p\nin \bigcup_{\al<\omi}W_\al$. 
Therefore 
it follows from (*) that $p\in\sid Xt$ is 
a $\ip{}1$ relation over $\hc$. 
This proves the lemma.
}%
\epf

\parf{The model}
\las{bgm}

Here we start the proof of Theorem~\ref{mt}.
The key idea of \cite[Part C]{h74} consists in 
making use of the \dd\omi long iterated extension 
of $\rL$, where the forcing at each step is the 
finite-support product of all elements of the 
fan $\fan\gQ$ defined within the extension 
obtained at the previous step of the iteration. 
We are going to define such an extension as a 
submodel of a more elementary  
background set universe $\gM$.

To define the latter, we 
consider the forcing notion $\gQ^{\omi\ti\omi}\in\rL$ 
(finite support). 
As $\gQ=Q^\om$, the forcing $\gQ^{\omi\ti\omi}$ is 
order-isomorphic to $Q^{\omi},$ of course.
The forcing $\gQ^{\omi\ti\omi}\in\rL$ 
naturally adjoins an array
%$\sis{g_{\nu\ga}}{\nu,\ga<\omi}$ 
of mutually 
\dd\gQ generic sets $G^{\nu\ga}\sq\gQ$, 
$\nu,\ga<\omi$ 
%(as in \ref{gq}) 
to $\rL$. 
%, such that 
%$T(g_{\nu\ga},a_{\nu\ga})$ holds for all $\nu,\ga$.
We let $\gM$, \rit{the background model}, be the extension 
$\rL[\sis{G^{\nu\ga}}{\nu,\ga<\omi}]$.

If $u\in\rL$, $u\sq\omi\ti\omi$, then let 
$\gM\ret u=\rL[\sis{G^{\mu\da}}{\ang{\mu,\da}\in u}]$.
In particular, if $\nu,\ga<\omi$ then put 
$\gM_\nu=\gM\ret{\nu\ti\omi}
=\rL[\sis{G^{\mu\da}}{\mu<\nu,\da<\omi}]$ 
and $\gM_{\nu\ga}=\gM\ret{\nu\ti\ga}$.

\ble
[by \ref{gq}(C)]
\lam{omi}
$\gM$ preserves all\/ \dd\rL cardinals. 
If\/ $x\in\gM$ is a real then\/ $x$ belongs to 
some\/ $\gM_{\nu\ga}$, $\nu\yi\ga<\omi$.  
Every\/ $\gM_{\nu\ga}$ is a\/ \dd\gQ generic 
extension of\/ $\rL$.\qed
\ele
        
{\ubf The actual model for Theorem~\ref{mt} 
will be a certain 
subuniverse $\gN\sq\gM$.} 
%, which partially inherits the product structure.

\bdf
\lam{amo}
Arguing in\/ $\gM$, we define, by transfinite induction 
on $\nu$, an array of 
countable \dd\gQ names $\bt_{\nu\ga}\sq\gQ\ti\om$, 
such that 
\ben
\nenu
\itlb{amo2}%
if $\nu,\ga<\omi$ then 
$\sis{\bt_{\mu\da}}{\mu\le\nu,\da<\omi}\in\gM_\nu$, and   
$\sis{\bt_{\mu\da}}{\mu\le\nu,\da<\ga}\in\gM_{\nu\ga}$, 
so that each particular $\bt_{\nu\ga}$ belongs to 
$\gM_{\nu,\ga+1}$.
\een
We also define \rit{derived objects}, namely 
\ben
\nenu
\atc
\itlb{amo1}%
reals $\br_{\nu\ga}= \bt_{\nu\ga}[G^{\nu\ga}]\sq\om$, 
sets 
$\bsg_{\nu\ga}=\sid{\br_{\nu\ga}}{\bt_{\nu\ga}}\sq\gQ$, 
forcing notions 
$\gQ_{\nu\ga}=\stk\gQ{\tle{\bt_{\nu\ga}}}$, and 

\itlb{amo1+}% 
models \ 
$\gN_\nu=\rL[\sis{\bt_{\mu\da},\br_{\mu\da}}{\mu<\nu,\da<\omi}]$, 
\
$\gN_{\nu\ga}=\rL[\sis{\bt_{\mu\da},\br_{\mu\da}}{\mu<\nu,\da<\ga}]$  
\ ($\ga<\omi$);
\een  
which, by construction and the results of Section~\ref{inte},
%\dd{\gQ_{\nu\ga}}generic 
satisfy the following: 
\ben
\nenu
\atc\atc\atc
\itlb{amo4}%
$\br_{\nu\ga}\in\gM_{\nu+1,\ga+1}[G^{\nu\ga}]$, 
$\gQ_{\nu\ga}=\stk\gQ{\tle{\bt_{\nu\ga}}}$ is a forcing 
notion in $\gM_{\nu,\ga+1}$, and 
$\bsg_{\nu\ga}\sq\gQ$ is a set \dd{\gQ_{\nu\ga}}generic 
over $\gM_{\nu,\ga+1}$ and over $\gM_{\nu}$, and satisfying 
$\rL[\bt_{\nu\ga},\bsg_{\nu\ga}]=\rL[\bt_{\nu\ga},\br_{\nu\ga}]$,
by Corollary~\ref{ptsr}\;;

\itlb{amo5}%
if $\nu,\ga<\omi$ then 
the arrays  
$\sis{\tle{\bt_{\mu\da}},\gQ_{\mu\da}}{\mu\le\nu,\da<\omi}$, 
$\sis{\bsg_{\mu\da},\br_{\mu\da}}
{\mu<\nu,\da<\omi}$
belong to  
$\gM_\nu$, and the arrays  
$\sis{\tle{\bt_{\mu\da}},\gQ_{\mu\da}}{\mu\le\nu,\da<\ga}$, 
$\sis{\bsg_{\mu\da},\br_{\mu\da}}{\mu<\nu,\da<\ga}$
belong to  
$\gM_{\nu\ga}$;

\itlb{amo6}%
therefore $\gN_\nu\sq\gM_\nu$ and 
$\gN_{\nu\ga}\sq\gM_{\nu\ga}$, $\kaz\ga$.
\een

Now the step. 
Suppose that $\nu<\omi$ and all sets 
$\bt_{\mu\da}\sq\gQ\ti\om$ 
%,   
($\mu<\nu\yi\da<\omi$)
are defined, so that 
$\sis{\bt_{\mu\da}}{\mu<\nu,\da<\omi}\in\gM_\nu$, 
and if $\ga<\omi$ then 
$\sis{\bt_{\mu\da}}{\mu<\nu,\da<\ga}\in\gM_{\nu\ga}$; 
this is slightly weaker than \ref{amo2} since does 
not include $\nu$ itself. 
Then $\gQ_{\mu\ga}$, $\bsg_{\mu\ga}$, $\br_{\mu\ga}$, 
$\gN_\nu$, $\gN_{\nu\ga}$ 
as in \ref{amo1}, \ref{amo1+} are defined as well. 
The goal is to define $\bt_{\nu\ga}\yi \ga<\omi$.

Note that 
$\gM_\nu=\rL[\sis{G^{\mu\ga}}{\mu<\nu,\ga<\omi}]$ 
is a \dd{\gQ^{\omi}}generic extension of $\rL$, 
hence GCH is true in $\gM_\nu$, 
and hence in $\gN_\nu\sq\gM_\nu$ as well. 
Therefore it holds in 
%$\gN_\nu$ and any 
$\gN_{\nu\xi}$ 
that there exist only \dd{\ali}many countable 
sets $t\sq\gQ\ti\om$; 
let $\sis{t^{\xi}_{\nu\eta}}{\eta<\omi}$ be 
the G\"odel-least 
(relative to 
%the array 
$\sis{\bt_{\mu\da},\br_{\mu\da}}{\mu<\nu,\da<\xi}$ 
as the parameter)
enumeration of all such $t$ in $\gN_{\nu\xi}$.

Let $\Om=\ens{\ga+1}{\ga<\omi}$ 
(successor ordinals). 
Fix a bijection  
$\bi:\Om\text{ onto }\omi\ti\omi$, $\bi\in\rL$, 
satisfying $\bi\obr(\xi,\eta)>\max\ans{\xi,\eta}$ 
for all $\xi,\eta<\omi$. 
If $\ga=\bi\obr(\xi,\eta)\in\Om$ then 
let $\bt_{\nu\ga}=t^{\xi}_{\nu\eta}$. 
Put $\bt_{\nu\ga}=\pu$ 
for all limit $\ga<\omi$.
The enumeration $\sis{\bt_{\nu\ga}}{\ga<\omi}$   
involves all at most countable sets $t\in\gN_\nu$, 
$t\sq\gQ\ti\om$, 
the whole sequence $\sis{\bt_{\nu\ga}}{\ga<\omi}$ belongs 
to $\gN_\nu$, and if $\ga<\omi$ then the subsequence   
$\sis{\bt_{\nu\da}}{\da<\ga}$ belongs to $\gN_{\nu\ga}$.
 
%If $\ga<\omi$ then we put $\bt_{\nu\ga}=t_{\nu\ga}$, and 
This ends the inductive step.

After the inductive construction is accomplished, 
we let 
$$
%\bay{lclcll}
%\gN &=&
%\rL[\sis{\bt_{\nu\ga},\br_{\nu\ga}}{\nu<\omi,\ga<\omi}] &=&
%\rL[\sis{\bt_{\nu\ga},\bsg_{\nu\ga}}{\nu<\omi,\ga<\omi}]\,,\\[1ex]
%
\gN\ret u 
=
\rL[\sis{\bt_{\nu\ga},\br_{\nu\ga}}{\ang{\nu,\ga}\in u}] 
=
\rL[\sis{\bt_{\nu\ga},\bsg_{\nu\ga}}{\ang{\nu,\ga}\in u}]\,,  
%\hspace*{-2ex}
\ \text{ for all $u\sq\omi\ti\omi$}\,,
%\eay
$$
and then $\gN=\gN\ret{\omi\ti\omi}$, 
and, equivalently to \ref{amo1+}, 
$\gN_\nu=\gN\ret{\nu\ti\omi}$ and  
$\gN_{\nu\ga}=\gN\ret{\nu\ti\ga}$  
for $\nu,\ga<\omi$.
We have by construction:
\ben
\nenu
\atc
\atc
\atc
\atc
\atc
\atc
\itlb{amo3}%
the whole sequence $\sis{\bt_{\nu\ga}}{\ga<\omi}$ belongs 
to $\gN_\nu$, and if $\ga<\omi$ then the subsequence   
$\sis{\bt_{\nu\da}}{\da<\ga}$ belongs to $\gN_{\nu\ga}$.
\qed
\een
\eDf

The next lemma explains further details.

\ble
\lam{maj}
Assume that\/ $\nu<\omi$, $t\in\gN_\nu$, 
$t\sq\gQ\ti\om$ is at most countable. 
%so\/ $\gQ_t=\stk\gQ{\tle t}\in\fan\gQ$ in\/ $\gN_\nu$. 

Then there is an ordinal\/ $\ga<\omi$ such that\/ 
$t=\bt_{\nu\ga}$.

In this case,\/ $\gQ_t=\gQ_{\nu\ga}$,  
the set\/ $\bsg_{\nu\ga}\in\gN_{\nu+1}$ 
is\/ \dd{\gQ_t}generic over\/ $\gM_\nu$, 
hence over the model\/ $\gN_\nu\sq\gM_\nu$ as well, 
and the real\/ 
$\br_{\nu\ga}= t[G^{\nu\ga}]$
satisfies\/ 
$\rL[\bt_{\nu\ga},\br_{\nu\ga}]=
\rL[\bt_{\nu\ga},\bsg_{\nu\ga}]$ 
and belongs to\/ $\gN_{\nu+1}$.
\ele
\bpf
Recall that $G^{\nu\ga}$ is \dd\gQ generic 
over $\gM_\nu$, hence over $\gN_\nu\sq\gM_\nu$ as well. 
It remains to use Theorem~\ref{32t}.
\epf

\bre
\lam{qq}
If $\nu<\omi$ then by construction the collection of 
all forcing notions $\gQ_{\nu\ga}$, $\ga<\omi$, 
is equal to the Harrington fan $\fan\gQ$ 
computed in $\gN_\nu$. 
Thus $\gN$ can be viewed as the \dd\omi long 
iterated \dd{\fan\gQ}generic extension of $\rL$. 

In particular, by \ref{qqp}(B), there is an index $\ga<\omi$ 
such that $\gQ_{\nu\ga}=\gQ$, hence, by Lemma~\ref{maj}, 
$\gN_{\nu+1}$ contains a set \dd\gQ generic over $\gM_\nu$ 
and over $\gN_\nu\sq\gM_\nu$. 
Similarly by \ref{qqp}(C) there is an index $\ga<\omi$ 
such that $\bt_{\nu\ga}=\bc\in\rL$ 
(see \ref{qqp}(C) on $\bc$),
and hence $\gQ_{\nu\ga}=\gQ_\bc$ adds a Cohen real 
$\br_{\nu\ga}\sq\om$ over 
$\gN_\nu$.
\ere

\parf{Key lemmas}
\las{keyl}

\vyk{
\ble
[in $\gN$]
\lam{kel0}
Let\/ $\nu<\omi$, $P\in\gN_\nu$, and it holds in\/ 
$\gN_\nu$ that\/ $P$ belongs to\/ $\fan\gQ$. 
 
\dd\gQ generic over\/ $\rL[x]$, hence there 
are reals\/ $g,a\in\gN$ such that the pair\/ 
$\ang{g,a}$ is\/ \dd Qgeneric over\/ $\rL[x]$.
\ele
\bpf
By Lemma~\ref{omi}, $x$ belongs to some 
$\gN_\nu\sq\gM_\nu$, $\nu<\omi$. 
By \ref{qq}, the submodel $\gN_{\nu+1}$, hence, 
$\gN$ as well, contains $G^{\nu\ga},$ 
a \dd\gQ generic set over $\gN_\nu$.
\epf
}

\ble
\lam{kel1}
If\/ $x\in\dn\cap\gN$, then there is a set\/ 
$G\in\gN$, 
\dd\gQ generic over\/ $\rL[x]$, hence there 
are reals\/ $g,a\in\gN$ such that the pair\/ 
$\ang{g,a}$ is\/ \dd Qgeneric over\/ $\rL[x]$.
\ele
\bpf
By Lemma~\ref{omi}, $x$ belongs to some 
$\gN_\nu\sq\gM_\nu$, $\nu<\omi$. 
By \ref{qq}, the submodel $\gN_{\nu+1}$ contains a set  
\dd\gQ generic set over $\gN_\nu$, 
hence over $\rL[x]$ as well.
\epf

\ble
\lam{kel1+}
Assume that\/ $\nu\yi\ga<\omi$, $u\in\rL$, 
$u\sneq\nu\ti\ga$, and\/ $W=(\nu\ti\ga)\bez u$. 
Then\/ 
$\gK=(\gN\ret u)[\sis{G^{\mu\da}}{\ang{\mu,\da}\in W}]$
is a\/ \dd\gQ generic extension of\/ $\gN\ret u$, 
and\/ $\gN_{\nu\ga}\sq\gK$.
\ele
\bpf
Let\/ $\gQ(\mu,\da)=\gQ$, 
for all\/ $\ang{\mu,\da}\in W$.
Note that $\gK$ is a 
\dd{\prod_{\ang{\mu,\da}\in W}\gQ(\mu,\da)}generic 
extension of $\gN\ret u$  
by construction, hence essentially a \dd{\gQ^\om}generic 
extension, yet $\gQ^\om$ is isomorphic to $\gQ$ 
as a forcing.
To prove $\gN_{\nu\ga}\sq\gK$, we check, by induction, 
that $\sis{\bt_{\ka\da},\br_{\ka\da}}{\da<\ga}\in\gK$ 
for all $\ka<\nu$. 
The induction hypothesis is $\ka<\nu$ and   
$\gN_{\ka\ga}\sq\gK$, and the goal is to 
``effectively'' prove that then 
$\sis{\bt_{\ka\da},\br_{\ka\da}}{\da<\ga}\in\gK$. 
We first remind that 
$\sis{\bt_{\ka\da}}{\da<\ga}\in\gN_{\ka\ga}$ 
by \ref{amo3} of Definition~\ref{amo}. 
Now, for any particular $\da<\ga$, if 
$\ang{\ka,\da}\in u$ then $\br_{\ka\da}$ belongs 
to $\gN\ret u$, hence, to $\gK$ as well, while if 
$\ang{\ka,\da}\in W$ then $G^{\ka\da}$ belongs 
to $\gK$, hence 
$\br_{\ka\da}=\bt_{\ka\da}[G^{\ka\da}]\in\gK$, 
as required. 
\epf

\bdf
[autonomous sets]
\lam{aur}
A set $x\in\gN\yi x\sq\rL$, is \rit{autonomous} if there is 
a countable set $u\in\rL\yi u\sq\omi\ti\omi$ such that 
$\rL[x]=\gN\ret u$.
\edf

\ble
[in $\gN$]
\lam{kelE}
If\/ $z$ is an \aut\ real and\/ $t\in\rL[z]$, 
$t\sq\gQ\ti\rL[z]$,
then there is a real\/ $b$ such that\/ 
$\ang{z,b}$ is \aut\ and\/ $\rL[z,b]=\rL[z][\Sg]$, 
where\/ $b=t[\Sg]$ and\/ $\Sg\sq\gQ$, $\Sg$ is\/ 
\dd{\gQ_t}generic over\/ $\rL[z]$.
\ele
\bpf
Let a countable $u\in\rL\yi u\sq\omi\ti\omi$ witness 
that $z$ is \aut. 
Then $u\sq\nu\ti\omi$ for some $\nu<\omi$, 
and $t\in\gN_\nu$. 
By Lemma~\ref{maj}, $t=\bt_{\nu\ga}$ for some $\ga$, 
and then $b=\br_{\nu\ga}$ is as required. 
To see that $\ang{z,b}$ is \aut\ note that 
$\rL[z,b]=\gN\ret v$, where $v=u\cup\ans{\ang{\nu,\ga}}$.
\epf

\ble
[in $\gN$]
\lam{kel2}
Let\/ $x\in\dn$ and\/ $\vpi(\cdot)$ be a\/ $\is13$ 
formula. 
Then\/ 
\ben
\renu
\itlb{kel2a}%
if\/ $\gQ$ forces\/ $\vpi(\dox)$ over\/ $\rL[x]$ 
then\/ $\vpi(x)$ is true\/ {\rm(in $\gN$)\,;}\snos
{In this lemma, $\dox=\ans{\bon}\ti x\in\rL[x]$ 
is a \dd\gQ name of 
$x$ itself.} 

\itlb{kel2b}%
if\/ $x$ is \aut\ and\/ $\vpi(x)$ is true\/ 
{\rm(in $\gN$)}
then\/ $\gQ$ forces\/ $\vpi(\dox)$ over\/ $\rL[x]$.
\een
\ele
\bpf
\ref{kel2a}
holds by Lemma~\ref{kel1}, 
since the truth of $\is13$ formulas 
passes to bigger models by Shoenfield.
To prove \ref{kel2b}, let 
$\rL[x]=\gN\ret u=
\rL[\sis{\bt_{\nu\ga},\br_{\nu\ga}}{\ang{\nu,\ga}\in u}]$, 
where $u\in\rL$, $u\sq\omi\ti\omi$ is countable.
Let $\vpi(\cdot)$ be $\sus z\,\psi(z,\cdot)$, $\psi$ 
being $\ip12$.
Assume that $\vpi(x)$ is true in $\gN$. 
There is a real $z\in\gN$ such that 
$\psi(z,x)$ is true in $\gN$.
There is an ordinal $\mu<\omi$, such that 
$u\sneq\mu\ti\mu$ and 
$z\in \gN_{\mu\mu}=
\rL[\sis{\bt_{\nu\ga},\br_{\nu\ga}}{\nu,\ga<\mu}]$. 
Then 
$z\in\gK=(\gN\ret u)[\sis{G^{\mu\da}}{\ang{\mu,\da}\in W}]$
by Lemma~\ref{kel1+}, where 
$W=(\mu\ti\mu)\bez u$. 
And $\gK$ is a \dd\gQ generic extension of 
$\gN\ret u=\rL[x]$ still by Lemma~\ref{kel1+}. 

On the other hand,  $\psi(z,x)$ is true in $\gK$ by 
Shoenfield, hence $\vpi(x)$ is true in $\gK$ as well. 
It follows that a condition in $\gQ$ 
forces $\vpi(\dox)$ over $\rL[x]$.
We conclude by Lemma~\ref{homle} that 
$\gQ$ forces\/ $\vpi(\dox)$ over\/ $\rL[x]$.
\epf

Recall that $\hc$ = hereditarily countable sets 
(in $\gN$).

\ble
\lam{kel6a}
If\/ $\vpi$ is a $\is{}n$ formula containing\/  
\dd\gQ names in $\HC$ 
{\rm(names of sets in $\HC[\uG]$)}, 
and if\/ $p\in\gQ$, then\/ 
$p\forn\gQ\HC\vpi$ 
is a\/ $\is\HC n$ assertion about\/ 
$p\yi\vpi$.
\ele
\bpf[sketch]
If $\vpi$ is a $\is12$ formula then, by the 
Mostowski absoluteness, $p\forn\gQ\HC\vpi$ 
iff $p\forn\gQ{M}\vpi$ over some countable 
transitive model $M$ of a sufficient fragment 
of $\zfc$, which is a $\is\hc1$ 
relation.
But, $\is\hc1$ relations are the same as $\is12$. 
This covers the case $n=1$ of the lemma. 

Step $\is{}n\to\ip{}n$. 
Suppose that 
$\vpi$ is a $\is{}n$ formula. 
Then $p\forn\gQ\HC\neg\:\vpi$ iff 
$\kaz q\,\big(q\leq p\imp\neg\:q\forn\gQ\HC\vpi\big)$. 
This leads to a $\ip{}n$ formula. 

Step $\ip{}n\to\is{}{n+1}$. 
Let 
$\vpi(x)$ be a $\ip{}n$ formula. 
Then $p\forn\gQ\HC\sus x\,\vpi(x)$ iff 
there is a \dd\gQ name $t\in\HC$ such that 
$p\forn\gQ\HC\vpi(t)$. 
This leads to a $\is{}{n+1}$ formula. 
\epf

\ble
\lam{kel6}
There is a recursive correspondence\/ $\vpi\mto\vpa$ 
between\/ $\is13$ formulas\/ 
{\rm(}and hence between\/ $\ip13$ formulas as well\/{\rm)} 
such that for all reals\/ $b\in\gN$, 
$\gQ$ forces\/ $\vpi(\dob)$ 
over\/ $\rL[b]$ if and only if\/ 
$\rL[b]\mo\vpa(b)$.
\ele
\bpf
Let $\jsg(z,t,p)$ be a $\ip{}1$ formula 
%that canonically  defines the set $\bbsg$ in $\lomi$ 
provided by Lemma~\ref{Sdef}. 

Given $\vpi$ a $\is13$ formula, define $\vpa(b)$ iff: 
$$
\TS
\sus t\in\rL\,\sus p\,\big(\,
\rL[b]\mo\jsg(b,t,p)\;\text{ and }\;
{p\form\gQ \vpi(\dob)}\;
\text{ over }\;\rL[b] 
\,\big).\snos
{$\dob=\ans{\bon}\ti b\in\rL[b]$ is the canonical 
\dd\gQ name of a real $b\sq\om$, where $\bon$ is the 
%canonically 
largest condition in $\gQ$.} 
$$
Prove that $\vpa$ is as required. 
Since $\is1n$ formulas correspond to the 
$\is{}{n-1}$ definability in $\HC$, $\vpa$ is 
a $\is\hc2$ formula by Lemma~\ref{kel6a}, 
%applied in $\rL[b]$, 
hence essentially a $\is13$ formula. 

Now suppose that $b$ is a real in $\gN$ and $\gQ$ forces 
$\vpi(\dob)$ over\/ $\rL[b]$. 
We have to prove that $\vpa(b)$ holds in $\rL[b]$. 
Note that $b\in\gN\sq\gM$. 
Hence $b\in\gM_{\nu\nu}$ for some $\nu<\omi$. 
But $\gM_{\nu\nu}$ is a \dd{\gQ^{\nu\ti\nu}}generic 
extension, hence, a \dd\gQ generic extension of $\rL$. 
Let say $\gM_{\nu\nu}=\rL[H]$, 
where $H\sq\gQ$, $H\in\gM$ is \dd\gQ generic over $\rL$; 
$H$ need not be in $\gN$.
Thus $b=t[H]$, where $t\in\rL\yi t\sq\gQ\ti\om$. 
Let $\Sg=\sid bt$; 
then $\Sg\sq\gQ$ and $\rL[\Sg]=\rL[b]$. 
Consider any $p\in\Sg$; thus $\jsg(b,t,p)$ holds in $\rL$. 
Under our assumptions, 
$p\form\gQ\vpi(\dob)$
over $\rL[b]$, hence we have $\vpa(b)$ in $\rL[b]$.

To prove the converse, assume that $\vpa(b)$ holds 
in $\rL[b]$, and this is witnessed by $t\in\rL$ and $p$. 
In particular $\jsg(b,t,p)$ holds, thus 
$p\in\Sg=\sid bt\sq \gQ$.  
Moreover, $p\form\gQ\vpi(\dob)$ over $\rL[b]$. 
It follows that if $G\sq\gQ$ is generic over $\rL[b]$ 
and $p\in G$ 
then $\vpi(b)$ is true in $\rL[b][G]$. 
Thus $\gQ$ forces $\vpi(\dob)$ 
over $\rL[b]$ by Lemma~\ref{homle}.
\epf
%{Lemma~\ref{kel6}}

%

\vyk{

\bcl
\lam{cl7}
For a term\/ $t$ which denotes a real in\/ $\rL^\gQ$, 
and for\/ $p\in\gQ$, let\/ $\vt(t,p,b)$ assert$:$ 
there is\/ $G$, an\/ \dd\rL generic filter on\/ $\gQ$, 
such that\/ $p\in G$ and\/ $b=t(G)$. 
Then\/ $\vt(t,p,b)$ is\/ $\ip\HC1$.
\ecl
\bpf[claim]
Let $\psi(t,p,b)\equiv \neg\:\vt(t,p,b)$. 
Using the recursion theorem, 
$\psi$ can be put into a $\is{}1$ form 
as follows: 
$\psi(t,p,b)$ iff: 
\ben
\tenu{(\arabic{enumi})}
\itla{z1}
there is $D\in W$ such that $\kaz q\in D$, 
$q$ compatable with $p$ 
$\imp$ $\kaz n\in\om\:(q\for t(n)\ne b(n))$, \ or

\itla{z2}
$\sus D\in W$ such that $\kaz q\in D$, $q$ compatable 
with $p$ 
$\imp$ there is an ordinal $\xi$ such that 
$\rL_\xi[b]\mo\psi(t,q,b)$ 
but $\rL_\xi[b]\not\mo\psi(t,p,b)$.
\een
The above just corresponds to the usual inductive 
definition of 
those $p$ which are inconsistent with the 
interpretation of $t$
as $b$.
\epF{Claim}

\begin{quote}
there is $t$ a term denoting a real in $\rL^\gQ$, 
there is $p\in\gQ$ such that $\vt(t,p,b)$, and such that 
$p\for\text{\lap{$\rL[t]^\gQ\mo\vpi(t)$}}$ 
(viewing $p$ as a member of $\gQ\ti\gQ$). 
\end{quote}

}

\parf{Reduction fails}
\las{nored}

In the remainder, we are going to prove that $\gN$ is a 
model for Theorem~\ref{mt}. 
The following is the first part of the proof.

\bte
\lam{tred}
In\/ $\gN$, there is a pair of\/ $\is13$ sets 
of reals, not reducible to a pair of boldface\/ 
$\fs13$ sets.
\ete
\bpf\snos
{Harrington relates the idea of the proof to Sami.} 
{\ubf Arguing in $\gN$}, consider the $\is13$ set 
$A=\ens{g}{\sus a\,T(g,a)}$,
and let $U\sq\dn$ be a $\is13$ set, 
universal in the sense that 
$\ens{e<\om}{e\we z\in U}\nin\ip13(z)$ 
for every $z\in\bn,$ where $e\we z$ adds  
$e\in\om$ as the leftmost term to $z\in\bn.$ 

Consider the $\is13$ sets $\om^\om\ti A$, $U\ti\om^\om$. 

Suppose to the contrary that, in $\gN$, 
there are $\fs13$ sets $A',U'\sq\dn,$ such that 
$$
A'\sq\om^\om\ti A,\;
U'\sq U\ti\om^\om,\;
A'\cap U'=\pu,\;
A'\cup U'=(\om^\om\ti A)\cup (U\ti\om^\om).
$$
Let $z$ be an \aut\ real such that $A',U'$ are $\is13(z)$.

\ble
[in $\gN$]
\lam{duy}
Assume that\/ $d\in\dn\bez U$. 
Then\/ $\ens{y}{\ang{d,y}\in A'}=A$ and there is\/ 
$y\in\dn$ such that\/ $\ang{d,y}\in A'$ and\/ 
$y$ is Cohen generic over $\rL[z,d]$.
\ele 
\bpf
If $\ang{d,y}\in A'$ then $y\in A$ by construction. 
Conversely assume that $y\in A$. 
Then $\ang{d,y}\in \bn\ti A$, but $\ang{d,y}\nin U\ti \bn$ 
(as $d\nin U$). 
Therefore $\ang{d,y}\in A'$ as required. 
\vyk{
$$ 
d\in\dn\bez U
\;\imp\;
\sus y\:(\ang{d,y}\in A' \,\land\,
\text{$y$ is \dd Cgeneric over $\rL[z]$})\,.
\eqno(\ast)
$$
}%
Prove the second claim. 
By Lemma~\ref{kel1}, 
there is a pair $\ang{g,a}$ in $\gN$, 
\dd Qgeneric over $\rL[z,d]$. 
Thus $g\in\dn$ is Cohen-generic over $\rL[z,d]$ 
while $a\sq\om$ satisfies $T(g,a)$, hence $g\in A$.
Thus we are done with $y=g$. 
\epF{Lemma~\ref{duy}}

\vyk{
Note that the right-hand side of $(\ast)$ is $\is13(z,d)$. 
Since $\ens{\ang{z,e}}{\ang{z,e}\notin U\land e\in\om}$ is not 
$\is13(z)$, we have that there is an integer $e$ such that 
$$ 
\sus y\:(\ang{\ang{z,e},y}\in A' \,\land\,
\text{$y$ is \dd Cgeneric over $\rL[z]$})\,,
$$
and such that $\ang{z,e}\in U$.
}

\vyk{
\ble
[in $\gN$]
\lam{due}
Assume that\/ $d\in\dn\bez U$. 
Then\/ $\ens{y}{\ang{d,y}\in A'}=A$ and there is\/ 
$y\in\dn$ such that\/ $\ang{d,y}\in A'$ and\/ 
$y$ is \dd Cgeneric over $\rL[z,d]$.
\ele 
\bpf

\epf
}

Consider the sets
$K=\ens{e\we z}{e<\om \land e\we z\nin U}\nin\is13(z)$ 
and   
$$
K'=\ens{e\we z}
{\sus y\:(e<\om \land\:\ang{\ang{z,e},y}\in A' \,\land\,
\text{$y$ is \dd Cgeneric over $\rL[z]$})}.
$$
Clearly $K\nin\is13(z)$ by the choice of $U$, 
while $K'$ is $\is13(z)$, and we have $K\sq K'$ by 
Lemma~\ref{duy}. 
Thus $K\sneq K'$. 
We conclude that there is an integer $e$ such that 
$e\we z\in U$ and $e\we z\in K'$, so that 
$$ 
\sus y\:(\ang{e\we z,y}\in A' \,\land\,
\text{$y$ is Cohen-generic over $\rL[z]$})\,.
$$
Fix such a number $e<\om$.%

Let 
$A''=\ens{y}{\ang{e\we z,y}\in A'}$; $A''\sq A$. 
We claim that 
$A''$ is $\id13(z)$. 
%(since $y\nin A''\eqv \ang{\ang{z,e},y}\in U'$).
Indeed $y\nin A''\eqv \ang{e\we z,y}\nin A'$. 
But $e\we z\in U$, hence $\ang{e\we z,y}\in U\ti\bn.$ 
Therefore 
$\ang{e\we z,y}\nin A'\eqv\ang{e\we z,y}\in U'$. 
This yields the 
%(less trivial) 
$\ip13$ definition for $A''$.

Let $\vpi$ be a $\ip13$ formula such that 
$y\in A''\eqv \vpi(z,y)$ in $\gN$. 
%Notice that $A''\sq A$. 
By the choice of $e$, there is a real 
$g\in A''$, Cohen-generic over $\rL[z]$. 
So $\vpi(z,g)$ is true in $\gN$.
It follows by Lemma~\ref{kel2} that  
$\gQ$ forces\/ $\vpi(\dotz,\dog)$ over\/ $\rL[z,g]$.
%$\rL[z,g]^\gQ\mo\vpi(z,g)$.

So by Lemma~\ref{kel6} we have $\rL[z,g]\mo\vpa(z,g)$.
By the genericity of $g$, there is 
a Cohen condition $\bag\in \dC\yt \bag\sq g$ such that 
$\bag\for
\text{\lap{
$\rL[z,\dog]\mo\vpa(z,\dog)$
}}
$. 

Recall that $z$ is autonomous. 
Let this be witnessed by a countable   
$u\in\rL\yi u\sq\nu\ti\vt$, where $\nu,\vt<\omi$; 
thus 
$\rL[z]=\gN\ret u=
\rL[\sis{\bt_{\mu\da},\br_{\mu\da}}{\ang{\mu,\da}\in u}]
\sq\gM_{\nu\vt}\sq\gM_\nu$. 
By \ref{qq}, there is an ordinal $\ga<\omi$ such 
that $\bt_{\nu\ga}=\bc$ and 
$\gQ_{\nu\ga}$ adds a Cohen real 
over $\gM_\nu$, so $\br_{\nu\ga}\in\dn$ is 
a Cohen real over $\gM_\nu$. 
Changing appropriately a finite number of values 
$\br_{\nu\ga}(k)$, we get another real $g'\in\dn,$ 
Cohen-generic over $\gM_\nu$, and satisfying 
$\bag\sq g'$ and still 
$(\gN\ret u)[g']=(\gN\ret u)[\br_{\nu\ga}]$, 
or in other words, 
$\rL[z,g']= (\gN\ret u)[g']= \gN\ret v$, 
where $v=u\cup\ans{\ang{\nu,\ga}}$.
Thus $\ang{z,g'}$ is \aut.

To conclude, $\ang{z,g'}$ is \aut\ and 
$\rL[z,g']\mo\vpa(z,g')$ by the choice of $g'$ 
and of $\bag\in\dC$.
It follows by Lemma~\ref{kel6} that $\vpi(z,g')$ is 
true in some/every \dd\gQ generic extension of $\rL[z,g']$. 
Therefore $\vpi(z,g')$ is true in $\gN$   
by Shoenfield.

%So by Claim~\ref{cl6} and Claim~\ref{cl5}, $M\mo\vpi(z,g')$.
Thus $g'\in A''$ by the choice of $\vpi$, 
hence $g'\in A$ and $\sus a\,T(g',a)$ holds in $\gN$.
Thus $\gQ$ forces $\sus a\,T(g',a)$ over $\rL[z,g']$ 
by Lemma~\ref{kel2}\ref{kel2b}, as $\ang{z,g'}$ is \aut. 

\ble
\lam{qfor}
$\gQ$ forces\/ $\sus a\,T(g',a)$ over\/ $\rL[g']$.
\ele
\bpf
It suffices to get a \dd\gQ generic 
extension of $\rL[g']$ in which $\sus a\,T(g',a)$ holds. 
Consider a set $G',$ \dd\gQ generic over $\rL[z,g']$. 
Then $\sus a\,T(g',a)$ holds in $\rL[z,g',G']$ by the 
above. 
Recall that $z$ belongs to a \dd\gQ generic extension 
$\gM_{\nu\vt}\sq\gM_\nu$ of $\rL$. 
Therefore $\rL[z]$ is a \dd{\gQ_t}generic extension 
of $\rL$ by Lemma~\ref{sust}, where $t\sq\gQ\ti\om$ 
is countable.
In other words, $\rL[z]=\rL[\Sg]$, where $\Sg\sq\gQ$ 
is \dd{\gQ_t}generic over $\rL$.

On the other hand, $g'$ is Cohen-generic over $\rL[z]$ 
while $G'$ is \dd\gQ generic over $\rL[z,g']$. 
It follows, by the product forcing theorem, that 
$G'$ is \dd\gQ generic over $\rL[g']$ and 
$\rL[z,g',G']=\rL[g',G'][\Sg]$ is a \dd{\gQ_t}generic 
extension of $\rL[g',G']$. 

Let $G\sq\Sg$ be \dd\Sg generic over $\rL[g',G'][\Sg]$. 
Then $\rL[g',G',\Sg,G]$ is a \dd\gQ generic 
extension of $\rL[g',G']$ by Theorem~\ref{32t}\ref{32t3}, 
hence a \dd\gQ generic extension of $\rL[g']$ as well 
because $\gQ\ti\gQ$ is order-isomorphic to $\gQ$.
Finally $\sus a\,T(g',a)$ holds in $\rL[g',G',\Sg,G]$ 
by Shoenfield as it holds in $\rL[z,g',G']$, 
a smaller model.
\epF{Lemma} 

But this contradicts Lemma~\ref{c1}.
\epF{Theorem~\ref{tred}}

\parf{Separation holds}
\las{+sep}

The next theorem is the final part of the proof of 
Theorem~\ref{mt}.

\bte
\lam{tsep}
In\/ $\gN$, any pair of disjoint\/ $\fp13$ sets 
is separable by a\/ $\fd13$ set.
\ete
\bpf
\rit{We argue in\/ $\gN$.}
Let $\vpi_0(\cdot,\cdot)\yd\vpi_1(\cdot,\cdot)$ 
be $\ip13$ formulas and $z$ be an autonomous real 
parameter, such that   
the $\fp13$ sets 
$A_i=\ens{x}{\vpi_i(z,x)}$, $i=0,1$, 
%$A_1=\ens{x}{\vpi_1(z,x)}$ 
are disjoint. 
Let $\vpi_i^\ast(z,x)$ be the \dd{\ip13}formulas 
provided by Lemma~\ref{kel6}. 
% $\kaz x\,(\vpi_0(z,x)\imp\neg\:\vpi_1(z,x))$ 
Then the sets  
$B_i =
\ens{x}{\rL[z,x]\mo\vpi_i^\ast(z,x)}$, 
$i=0,1$, are $\fp13$, and   
$A_i\sq B_i$ by Lemma \ref{kel2}.

We claim that $B_0\cap B_1=\pu$. 
Assume to the contrary that $b\in B_0\cap B_1$. 
The goal of the following argument is to get another 
real $b'\in B_0\cap B_1$, with the extra property that 
$\ang{z,b'}$ is \aut.

As $z$ is \aut, Lemma~\ref{kel1+} implies that $b$ 
belongs to a \dd\gQ generic extension 
$\rL[z][G]$ of $\rL[z]$. 
Therefore, by Lemma~\ref{sust}, there is a countable 
set $t\sq\gQ\ti\om$, $t\in\rL[z]$, 
such that $b=t[G]$ and $\rL[z,b]$ is a 
\dd{\gQ_t}generic extension of $\rL[z]$.    
By Lemma~\ref{homle}\ref{homle2}, $\gQ$ forces 
\lap{$\vpi_0^\ast(\doz,t[\uG])$ and 
$\vpi_1^\ast(\doz,t[\uG])$ 
hold in $\rL[z,t[\uG]]$} over $\rL[z]$, therefore 
$\gQ_t$ forces 
\lap{$\vpi_0^\ast(\doz,t[\uG])\land \vpi_0^\ast(\doz,t[\uG])$} 
over $\rL[z]$. 
On the other hand, by Lemma~\ref{kelE}, 
there is a real $b'$ such that $\ang{z,b'}$ is \aut\  
and $\rL[z,b']=\rL[z][\Sg]$, where $b'=t[\Sg]$ and 
$\Sg\sq\gQ$ is \dd{\gQ_t}generic over $\rL[z]$. 
Thus 
%the sentences 
$\vpi_0^\ast(z,b')$ and $\vpi_1^\ast(z,b')$ 
hold in $\rL[z,b']$, 
hence $b'\in B_0\cap B_1$, and $\ang{z,b'}$ is \aut.

But then $b'\in A_0\cap A_1$ 
by Lemmas \ref{kel6} and \ref{kel2}, 
%$\gN\mo \vpi_0(z,b')\land\vpi_1(z,b')$.
contradiction. 

Thus we have $B_0\cap B_1=\pu$. 
Then we separate 
$B_0$ from $B_1$ by a $\id13(p)$ set of reals by 
a standard argument. 
Indeed let $\gle zx$ be the canonical \lap{good} 
G\"odel wellordering of the reals in $\rL[z,x]$.
Let $\sus y\,\vt_0(z,x,y)$ be the canonical transformation 
of $\neg\:\vpi_1^\ast(z,x)$ to \dd{\is13}form, and 
$\sus y\,\vt_1(z,x,y)$ be the canonical transformation 
of $\neg\:\vpi_0^\ast(z,x)$ to \dd{\is13}form, so that 
$\vt_i$ are \dd{\ip12}formulas and
\bce 
$B_i\sq C_i=\ens{x}{\rL[x,z]\mo \sus y\,\vt_i(z,x,y)}$, 
\ \ \ 
$i=0,1$. 
\ece
Here $C_i$ is the complement to $B_{1-i}$, thus 
$C_0\cup C_1=$ all reals. 
If $x\in C_i$ then let $y_i(x)$ be the \dd{\gle zx}least 
real $y$ satisfying $\rL[x,z]\mo \vt_i(z,x,y)$. 
The sets 
$$
\bay{rcl}
D_0
&=&\ens{x\in C_0}{x\nin C_1\lor y_0(x)\gle zxy_1(x)}\\[1ex]
D_1
&=&\ens{x\in C_1}{x\nin C_0\lor y_1(x)\gl zxy_0(x)} 
\eay
$$
then satisfy $A_i\sq B_i\sq D_i\sq C_i$ and 
$D_0\cup D_1=$ all reals. 
On the other hand, a standard argument 
(as in the proof of \dd{\fs1n}reduction in $\rL$) 
shows that both $D_i$ are $\fs13$ sets. 
It follows that $D_0$ is a $\fd13$ set separating 
$A_0$ from $A_1$. 
\epf

\parf{Comments and questions}
\las{comm}

We may note the following substantial 
inventions in Harrington's proof.
\bit
\item
The \rit{localization property} in $\gN$, that is, 
the reduction of the truth of a formula $\vpi(x)$ 
in the final model $\gN$, first, to the truth in 
\dd\gQ generic extensions of $\rL[x]$ by Lemma~\ref{kel2}, 
and second, to the truth in $\rL[x]$ itself by 
Lemma~\ref{kel6}. 
This is quite similar to the \lap{important lemma} 
of Solovay~\cite[page 18]{sol}, but achieved in a much less 
friendly generic model.

\item
The \rit{Harrington fan} construction of \ref{qqp}  
(see also Remark~\ref{qq}) 
which allows to inhibit (by Lemma~\ref{kelE}) 
the fact that Lemma~\ref{kel2}\ref{kel2b} holds 
only for autonomous reals. 

\item
Unlike the Separation counterexamples in specific 
models in \cite[Part B]{h74} or say \cite{kl28}, 
the Reduction counterexample as in Theorem~\ref{tred} 
(which Harrington grants to Sami) 
is not something explicitly designed by the intended 
definability structure of generic reals 
in the model considered. 
\eit
Harrington ends \cite[Part C]{h74} with the following 
remark:
\begin{quote}
We believe that this result 
[= Theorem~\ref{mt}] 
can be generalized by replacing $3$ by 
any integer $n\ge3$. 
We also believe that that this result can be improved so as 
to obtain a model of $\ZFC$ in which both $\sep(\fp13,\fd13)$ 
and $\sep(\fs13,\fd13)$ hold.
At the moment though these beliefs are just expressions of faith 
(or is it hope?).
\end{quote}
The second part of this \lap{expressions of faith or hope} 
was partially materialized in \cite[Part C]{h74}, where, 
for an arbitrary $n\ge3$,  
a model of\/ $\ZFC$ is presented, in which 
$\sep(\ip1n,\id1n)$ 
and $\sep(\is1n,\id1n)$ (note the lightface classes!) 
both hold for sets of integers.
(The proof is given for $n=3$ only.)
The rest presumably remains as open as it was in 1970s.
%Regarding the very recent results, 
%a model, in which Uniformization 
%(and hence Reduction) holds for $\fp13$, was discussed 
%in \cite{sdfun}.

%\np

\parf{Reduction holds in extensions by Cohen reals}
\las{redf}

Here we sketch 
{\ubf the proof of Claim \ref{mt'2} 
of Theorem \ref{mt'}.}
%Proofs of the other two claims will appear elsewhere.

Let the set universe $\rV$ be an extension 
of $\rL$ by a transfinite sequence of 
Cohen-generic reals. 
The following is a known property of 
the Cohen forcing $\dC=\bse$ 
and Cohen extensions.

\ble
\lam{kogen}
If\/ $x,y\in\rV$ are reals then either\/ $y\in\rL[x]$ 
or\/ $\rL[x,y]$ is a\/ \dd\dC generic extension of\/ 
$\rL[x]$. 
In particular, $y$ belongs to a\/ \dd\dC generic 
extension of\/ $\rL[x]$.\qed
\ele

{\ubf\boldmath Case $n=3$.} 
%It is known that if $x,y\in\rV$ are reals then 
%there is a real $z\in\rV$, Cohen-generic over 
%$\rL[x]$ and satisfying $y\in\rL[x,z]$. 
We claim 
%(by using Shoenfield and the homogeneity of the 
%Cohen forcing $\dC$) 
that if $\vpi(x)$ is a $\is13$ or $\ip13$ formula 
with $x$ as the only real parameter, then
\bce
\qquad$\vpi(x)$ holds in $\rV$ \ \ \  iff \ \ \  
$\rL[x]\mo\underbrace{\text{%
$\La$ \dd\dC forces $\vpi(\dox)$ 
over the universe}}_{\vpa(x)}$.
\hfill(1)
\ece 
Here $\La\in\dC$ (the empty sequence) 
is the weakest Cohen condition, and $\dox=\ans\La\ti x$ 
is the canonical Cohen name for a set $x$ in the 
ground model. 

To prove (1) for a $\is13$ formula 
$\vpi(x):=\sus x_1\psi(x,y)$, $\psi$ being $\ip12$, 
assume that $\rV\mo\vpi(x)$, 
hence there is a real 
$y\in\rV$ satisfying $\psi(x,y)$. 
But $y$ belongs to a \dd\dC generic extension of 
$\rL[x][g]$ of $\rL[x]$ by Lemma~\ref{kogen}. 
Then $\sus y\,\psi(x,y)$ is true in $\rL[x,g]$ 
by Shoenfield, and hence $\sus y\,\psi(\dox,y)$ 
is \dd\dC forced by $\La$ over $\rL[x]$ 
(by the homogeneity of $\dC$), 
that is, $\rL[x]\mo \vpa(x)$. 
Conversely let $\rL[x]\mo\vpa(x)$. 
Let $g\in\rV$ be a \dd\dC generic real over $\rL[x]$. 
Then $\rL[x,g]\mo\vpi(x)$, hence 
$\rV\mo\vpi(x)$ by Shoenfield. 

To check (1) for a $\ip13$ formula 
$\Phi(x):=\neg\,\vpi(x)$, 
$\vpi$ being $\is13$, assume first that 
$\rV\mo\Phi(x)$. 
Then $\rV\not\mo\vpi(x)$, hence by (1) 
$\vpi(\dox)$ is not \dd\dC forced 
by $\La$ over $\rL[x]$, thus by the homogeneity 
$\Phi(\dox)$ is forced, that is, 
$\rL[x]\mo \Pha(\dox)$.
Conversely if $\rL[x]\mo \Pha(\dox)$, then 
definitely $\rL[x]\not\mo \vpa(\dox)$, 
thus $\rV\not\mo \vpi(x)$, and hence 
$\rV\mo \Phi(x)$. 

Pretty similar to the proof 
of Lemma~\ref{kel6a}, $\vpa(\cdot)$ is a formula of 
type $\is13$, resp., $\ip13$ formula provided $\vpi$ 
itself is of this type. 
The next lemma will be used below.

\ble
\lam{xgh}
Let\/ $\vpi(x)$ be a\/ $\is13$ or\/ $\ip13$ formula. 
Let\/ $g$ be a real\/ \dd\dC generic over\/ $\rL[x]$, 
and\/ $\rL[x,g]\mo\vpa(x)$. 
Then\/ $\rL[x]\mo\vpa(x)$.
\ele
\bpf
By the homogeneity, it suffices to get a 
\dd\dC generic real\/ $h$ over\/ $\rL[x]$, such that 
$\rL[x,h]\mo \vpi(x)$. 
Let $g'$ be \dd\dC generic over\/ $\rL[x,g]$. 
Then $\rL[x,g,g']\mo \vpi(x)$, so it remains to 
make use of $h=\ang{g,g'}$.
\epf

Now, consider $\fs13$ sets 
$A_0=\ens{x}{\vpi_0(x,z)}$ and 
$A_1=\ens{x}{\vpi_1(x,z)}$ 
in $\rV$, where $\vpi_i$ are $\is13$ formulas. 
Then $A_i=\ens{x}{\rL[x,z]\mo\vpa_i(x,z)}$ by the above, 
where $\vpa_i(x,z):=\sus y\Phi_i(x,z,y)$ 
are $\is13$ formulas by the claim, so $\Phi_i$ are 
$\ip12$ formulas.

If $x\in A_i$ then let $y_i(x,z)$ be the \dd{\gle xz}least 
real $y$ satisfying $\rL[x,z]\mo \Phi_i(x,z,y)$. 
The sets 
$$
\bay{rcl}
B_0
&=&\ens{x\in A_0}{x\nin A_1\lor y_0(x,z)\gle xz y_1(x,z)}\\[1ex]
B_1
&=&\ens{x\in A_1}{x\nin A_0\lor y_1(x,z)\gl xz y_0(x,z)} 
\eay
$$
then satisfy $B_i\sq A_i$ and 
$B_0\cap B_1=\pu$, and belong to $\fs13$. 

{\ubf\boldmath Case $n\ge4$.} 
Let say $n=4$ exactly; it will be clear how to treat 
the general case. 
Suppose that $\vpi(x):=\sus y\,\psi(x,y)$ 
is a $\is14$ formula, 
$\psi$ being $\ip13$. 
Then a $\ip13$ formula $\psa(x,y)$ has been defined as 
above, such that $\rV\mo\psi(x,y)$ iff 
$\rL[x,y]\mo\psa(x,y)$, for all reals $x,y\in\rV$.
We claim that then
\bce
\qquad$\rV\mo\vpi(x)$  \ \ \  iff \ \ \  
$\rL[x]\mo\underbrace{%
\text{$\La$ \dd\dC forces \lap{$\sus y\,\psa(\dox,y)$} 
over the universe}}_{\vpa(x)}$.
\hfill(2)
\ece 
Indeed assume that $\rV\mo\vpi(x)$, 
hence there is a real 
$y\in\rV$ satisfying $\psi(x,y)$. 
It follows that $\rL[x,y]\mo\psa(x,y)$. 
But $y$ belongs to a \dd\dC generic extension of 
$\rL[x]$ by Lemma~\ref{kogen}. 
Therefore, $\sus y\,\psa(\dox,y)$ 
is \dd\dC forced by $\La$ over $\rL[x]$ 
(by the homogeneity of $\dC$), 
that is, $\rL[x]\mo \vpa(x)$. 

To prove the converse, let $\rL[x]\mo\vpa(x)$. 
Let $g\in\rV$ be a Cohen-generic real over $\rL[x]$. 
Then we have $\rL[x,g]\mo\sus y\,\psa(x,y)$. 
Let this be witnessed by a real $y\in \rL[x,g]$, 
thus $\rL[x,g]\mo\psa(x,y)$. 
However, by Lemma~\ref{kogen}, 
either $\rL[x,g]=\rL[x,y]$ --- and then 
$\rL[x,y]\mo\psa(x,y)$ and hence $\rV\mo\psi(x,y)$ 
and $\rV\mo\vpi(x)$, or $\rL[x,g]$ is a \dd\dC 
generic extension of $\rL[x,y]$ --- and then 
we still have $\rL[x,y]\mo\psa(x,y)$ by Lemma~\ref{xgh}, 
and then $\rV\mo\vpi(x)$ as just above. 

The proof of Reduction for a pair of \dd{\fs14}sets 
$A_0\yi A_1$ in $\rV$ goes on, on the base of (2), 
exactly as in the case $n=3$ above.

\vyk{
The principal ingredients are as follows. 

1. 
Let $\gM=\rL[\sis{x_\al}{\al<\omi}]$ is an extension 
of $\rL$ by \dd\ali many Cohen reals $x_\al$ then 
each submodel $\gM_\la=\rL[\sis{x_\al}{\al<\la}]$ 
($\la<\omi$) is a Cohen extension of $\rL$ 
(by a single Cohen real).

2. 
If a real $x$ is Cohen-generic over $\rL$ and 
$y\in\rL[x]$ is a real then 
\bit
\item
$y\in\rL$ or $\rL[y]$ is a Cohen extension of $\rL$;

\item
$x\in\rL[y]$ 
or $\rL[x]$ is a Cohen extension of $\rL[y]$.
\eit

It is well-known that if a real $x$ is Cohen-generic .
}

\bibliographystyle{plain}
{\small
%\bibliography{46,kle}
%

}

\end{document}